\pdfoutput=1
\RequirePackage{ifpdf}
\ifpdf 
\documentclass[pdftex]{sigma}
\else
\documentclass{sigma}
\fi

\numberwithin{equation}{section}

\newtheorem{Theorem}{Theorem}[section]
\newtheorem*{Theorem*}{Theorem}
\newtheorem{Corollary}[Theorem]{Corollary}
\newtheorem{Lemma}[Theorem]{Lemma}
\newtheorem{Proposition}[Theorem]{Proposition}
\newtheorem{Conjecture}[Theorem]{Conjecture}
 { \theoremstyle{definition}
\newtheorem{Definition}[Theorem]{Definition}

\newtheorem{Remark}[Theorem]{Remark} }

\newcommand{\ebar}{\overline{\varepsilon}}
\newcommand{\veps}{{\varepsilon}}
\newcommand{\NSL}{{\mathcal{L}}}
\newcommand{\BC}{{\mathbb{C}}}
\newcommand{\BR}{{\mathbb{R}}}
\newcommand{\BZ}{{\mathbb{Z}}}
\newcommand{\BN}{{\mathbb{N}}}

\newcommand{\wtJ}{{\widetilde{J}}}
\newcommand{\wtW}{{\widetilde{W}}}
\newcommand{\wtL}{{\widetilde{L}}}
\newcommand{\wtT}{{\widetilde{T}}}
\newcommand{\wtPi}{{\widetilde{\Pi}}}

\newcommand{\End}{{\textnormal{End}}}
\newcommand{\Spec}{{\textnormal{spec}}}

\begin{document}

\allowdisplaybreaks

\newcommand{\arXivNumber}{2211.01586}

\renewcommand{\PaperNumber}{063}

\FirstPageHeading

\ShortArticleName{Spectral Theory of the Nazarov--Sklyanin Lax Operator}

\ArticleName{Spectral Theory of the Nazarov--Sklyanin\\ Lax Operator}

\Author{Ryan MICKLER~$^{\rm a}$ and Alexander MOLL~$^{\rm b}$}

\AuthorNameForHeading{R.~Mickler and A.~Moll}

\Address{$^{\rm a)}$~Singulariti Research, Melbourne, Victoria, Australia}
\EmailD{\href{mailto:ry.mickler@gmail.com}{ry.mickler@gmail.com}}

\Address{$^{\rm b)}$~Department of Mathematics and Statistics, Reed College, Portland, Oregon, USA}
\EmailD{\href{mailto:amoll@reed.edu}{amoll@reed.edu}}
\URLaddressD{\url{https://alexander-moll.com}}

\ArticleDates{Received March 19, 2023, in final form August 27, 2023; Published online September 10, 2023}

\Abstract{In their study of Jack polynomials, Nazarov--Sklyanin introduced a remarkable new graded linear operator $\NSL\colon F[w] \rightarrow F[w]$ where $F$ is the ring of symmetric functions and~$w$ is a variable. In this paper, we (1) establish a cyclic decomposition $F[w] \cong \bigoplus_{\lambda} Z(j_{\lambda}, \NSL)$ into finite-dimensional $\NSL$-cyclic subspaces in which Jack polynomials $j_{\lambda}$ may be taken as cyclic vectors and (2) prove that the restriction of $\NSL$ to each $Z(j_{\lambda}, \NSL)$ has simple spectrum given by the anisotropic contents $[s]$ of the addable corners $s$ of the Young diagram of $\lambda$. Our proofs of (1) and (2) rely on the commutativity and spectral theorem for the integrable hierarchy associated to $\NSL$, both established by Nazarov--Sklyanin. Finally, we {conjecture that} the $\NSL$-eigenfunctions $\psi_{\lambda}^s {\in F[w]}$ {with eigenvalue $[s]$ and constant term} $\psi_{\lambda}^s|_{w=0} = j_{\lambda}$ are polynomials in the rescaled power sum basis $V_{\mu} w^l$ of $F[w]$ with integer coefficients.}


\Classification{05E05; 33D52; 37K10; 47B35}

\section{Introduction and statement of results}

Nazarov--Sklyanin introduced in~\cite{NaSk2} a graded linear operator $\NSL$ to the study of Jack polynomials~\cite{Jack, Mac, Stanley}. In Section~\ref{SUBSECNazarovSklyaninLaxOperator}, we recall the definition of $\NSL\colon F[w] \rightarrow F[w]$ in the polynomial ring~$F[w]$ where $F$ is the ring of symmetric functions with its usual grading and $\deg w=1$. They proved that if one considers the projection $\pi_0\colon F[w] \rightarrow F$ defined by setting $w=0$, then the operators $\mathcal{T}_{\ell} = \pi_0 \NSL^{\ell}$ pairwise commute in $F$ for all $\ell$ and are simultaneously diagonalized on Jack polynomials $j_{\lambda}$ with explicit eigenvalues. Moreover, as the second author observed in~\cite{Moll2, Moll5}, Nazarov--Sklyanin actually prove in~\cite{NaSk2} that the eigenvalues of $\mathcal{T}_{\ell}$ at $j_{\lambda}$ are precisely the $\ell^{\textnormal{th}}$ moments of the transition measure $\tau_{\lambda}$ of the anisotropic Young diagram of $\lambda$ studied in~\cite{Bi1, DoFe4, HoOb, Ke4, KvingeLicataMitchell, Las1, Ol1}.

 In this paper, we determine the spectrum of $\NSL$ in $F[w]$ and identify a distinguished polynomial basis $\psi_{\lambda}^s$ of eigenfunctions of $\NSL$ satisfying $\pi_0 \psi_{\lambda}^s=j_{\lambda}$. We state our result in Theorem~\ref{TheoremMAINRESULT} below. {In subsequent work~\cite{MicklerPAPER2}, the first author} uses the spectral theorem {established in the present paper to derive a new explicit system of constraints on Jack Littlewood--Richardson coefficients in terms of a simple new multiplication operation on partitions. Using the results of~\cite{MicklerPAPER2}, Alexandersson--Mickler~\cite{AlexanderssonMicklerPAPER3} prove new cases of the strong Stanley conjecture~\cite{Stanley}. We hope that the polynomials $\psi_{\lambda}^s$ introduced in this paper may inspire further applications and are also of independent interest. }

\subsection{The Nazarov--Sklyanin Lax operator} \label{SUBSECNazarovSklyaninLaxOperator} Consider the graded polynomial ring
\begin{equation} F = \BC[V_1, V_2, \ldots ] \label{H0Definition}\end{equation}
in which $\deg V_k = k.$ For any $\hbar \in \BC$, let $V_{-k}$ be the differential operator in $F$ defined by
\begin{equation} \label{VNegativeKDefinition} V_{-k} = \hbar k \frac{\partial}{\partial V_k}. \end{equation}

 Since $[V_{-k}, V_{k}] = \hbar k$,~(\ref{H0Definition}) is the Fock space representation of the Heisenberg algebra at level~$\hbar$. We will identify $F$ with the ring of symmetric functions in Section~\ref{SUBSECJackPolynomials}. Introduce a new variable $w$ with $\deg w =1$ and consider the graded polynomial and Laurent polynomial rings \begin{gather}\label{HDefinition}
 F[w] = \BC[w, V_1, V_2, V_3, \ldots], \\
 F\big[w, w^{-1}\big] = \BC\big[w, w^{-1}, V_1, V_2, V_3, \ldots\big]. \nonumber 
 \end{gather}
 Define the projection $\pi \colon F\big[ w, w^{-1}\big] \rightarrow F[w]$ to be the linear extension of
 \begin{equation}\label{SzegoProjectionFormula}
 \pi w^l = \begin{cases}
 \,\, w^l & \textnormal{if} \,\, l \geq 0, \\
 \,\, 0 & \textnormal{if} \,\, l < 0 .
 \end{cases}
 \end{equation}

 {We now present the Lax operator from~\cite{NaSk2} using the conventions from~\cite{Nak0, NekOk}.}

\begin{Definition} For any $\ebar, \hbar \in \BC$ and $\partial_w = \frac{\partial}{\partial w}$, the Nazarov--Sklyanin Lax operator $\NSL$ \textnormal{\cite{NaSk2}} is the graded linear operator $\NSL \colon F[w] \rightarrow F[w]$ in the graded polynomial ring $F[w]$ in {(\ref{HDefinition})} defined by
\begin{equation} \label{PresentationOfL} \NSL = \ebar w \partial_w + \sum_{k=1}^{\infty} V_{-k} w^k + \sum_{k=1}^{\infty} V_k\pi w^{-k} \end{equation}
for $V_{-k}$ in {(\ref{VNegativeKDefinition})}, $\pi$ in {(\ref{SzegoProjectionFormula})}. Below, we may write $\NSL = \NSL_{\ebar,\hbar}$ to emphasize that $V_{-k}$ depends on $\hbar$.\end{Definition}

Due to the presence of $\pi$ in~(\ref{PresentationOfL}), $\NSL$ is well defined with codomain $F[w]$. Moreover, $\NSL$ preserves total degree in $F[w]$ since $w \partial_w$ preserves degree, multiplication by $w^k$ and $V_k$ raise degree by $k$, and $V_{-k}$ and $\pi w^{-k}$ lower degree by $k$. The Lax operator $\NSL_{\ebar, \hbar}$ is a two-parameter perturbation of a~nilpotent linear operator $\NSL_{0,0}$ by linear differential operators which act as derivations in $F[w]$.

Since $\big\{w^{l}\big\}_{l=0}^{\infty}$ is a basis of $\BC[w]$ indexed by $\BN = \{0,1,2,\ldots\}$, $\NSL$ in~(\ref{PresentationOfL}) acts as an $\BN \times \BN$ matrix in $F[w] = F \otimes \BC[w]$ with coefficients in $\End(F)$: using the definition of $V_{-k}$ in~(\ref{VNegativeKDefinition}), this matrix is
\begin{equation} \label{MatrixPresentationOfL} \NSL =
\begin{bmatrix} 0 \ebar & V_1 & V_2 & V_3 & V_4 & \cdots \\
V_{-1} & 1 \ebar & V_1 & V_2 & V_3 & \ddots \\
V_{-2} & V_{-1} & 2 \ebar & V_1 & V_2 & \ddots \\
V_{-3} & V_{-2} & V_{-1} & 3 \ebar & V_1 & \ddots \\
V_{-4} & V_{-3} & V_{-2}& V_{-1} & 4 \ebar & \ddots \\
\vdots & \ddots & \ddots & \ddots & \ddots & \ddots
\end{bmatrix}.
\end{equation}

One can see from~(\ref{MatrixPresentationOfL}) that $\NSL_{0,0}$ is strictly upper triangular if $\ebar = \hbar =0$, hence nilpotent in~$F[w]$ as we already mentioned. If one sets $\ebar = 0$ but keeps $\hbar \neq 0$, the diagonal terms in~(\ref{MatrixPresentationOfL}) vanish and the entries of $\NSL_{0, \hbar}$ are constant along diagonals. In this special case, one may regard~$\NSL_{0, \hbar}$ as a block Toeplitz operator with infinite blocks from $\End(F)$.

\subsection{Jack symmetric functions} \label{SUBSECJackPolynomials} Recall that a partition $\lambda$ of $n$ is a sequence of non-negative integers $\lambda_i$ which are weakly decreasing $0 \leq \cdots \leq \lambda_3 \leq \lambda_2 \leq \lambda_1$ and satisfy $ \cdots + \lambda_3 + \lambda_2 + \lambda_1 =n$. {For any $\veps_1, \veps_2 \in \BC$, let $j_{\lambda} \in F$ be the Jack polynomial in the ring $F$ in~(\ref{H0Definition}) as defined in~\cite{LiQinWang, Nak0, NekOk, Okounkov2018ICM, QinBOOK}}. For example, the Jack polynomials in degrees $1 \leq |\lambda | \leq 3$ are
\begin{gather*} j_1 = V_1, \\
j_{2} = V_1^2 + \veps_1 V_2, \\
j_{1,1} = V_1^2 + \veps_2 V_2, \\
j_3 = V_1^3 + 3 \veps_1 V_1 V_2 + 2 \veps_1^2 V_3, \\
j_{1,2} = V_1^3 + (\veps_1 + \veps_2) V_1 V_2 + (\veps_1 \veps_2) V_3, \\
j_{1,1,1} = V_1^3 + 3 \veps_2 V_1 V_2 + 2 \veps_2^2 V_3.
\end{gather*}
For $\veps_2 \neq 0$, these polynomials are equivalent to $J_{\lambda}(p_1, p_2, \ldots ; \alpha)$, the integral form Jack symmetric function in~\cite{Jack, Mac, Stanley} where $p_k$ are the power sum symmetric functions and $\alpha$ is the Jack parameter. {Precisely, for $\veps_2 \neq 0$, if one sets} $p_k = V_k / (-\veps_2)$ and $\alpha = \veps_1 / (-\veps_2)$, then
\begin{equation} \label{JackPolynomialPresentation} j_{\lambda} = (- \veps_2)^{|\lambda|} J_{\lambda} \Bigg ( \frac{V_1}{(-\veps_2)}, \frac{V_2}{(-\veps_2)}, \ldots ; \frac{\veps_1}{ (- \veps_2)} \Bigg ).
\end{equation}
At $\veps_1 = \alpha$ and $\veps_2 = -1$, $j_{\lambda} = J_{\lambda}$.

 \subsection{Addable and removable corners of Young diagrams} 
 Any partition $\lambda$ determines a set
 \begin{equation} \label{DefinitionDLambda}
 \lambda= \bigcup_{r=1}^{\infty} \big\{ {(c-1, r-1)} \in \BN^2 \colon c \in \{1, 2, \ldots, \lambda_r\} \big \}
 \end{equation}
 called the Young diagram of $\lambda$. We use $\lambda$ to refer to either the sequence of parts $\lambda_i$ or to~(\ref{DefinitionDLambda}).
 \begin{Definition} 
 For $\lambda$ in {(\ref{DefinitionDLambda})}, define the addable and removable corner sets by
 \begin{gather} \label{AddSetFormula}
 \mathcal{A}_{\lambda} = \big\{ s \in \BN^2 \colon s \not\in \lambda \ \textnormal{and} \ \lambda \cup \{s\} \ \textnormal{is also a Young diagram} \big\}, \\
\label{RemoveSetFormula} \mathcal{R}_{\lambda} = \big\{ s \in \BN^2 \colon s \in \lambda \ \textnormal{and} \ \lambda \setminus \{s\} \ \textnormal{is also a Young diagram}\big\}. \end{gather}
It is also convenient to define the outer corner set $\mathcal{R}_{\lambda}^+ = \mathcal{R}_{\lambda} + (1,1)$ as shifts of removable corners. \end{Definition}
\begin{figure}[htb]

\centering
\includegraphics[width=0.42 \textwidth]{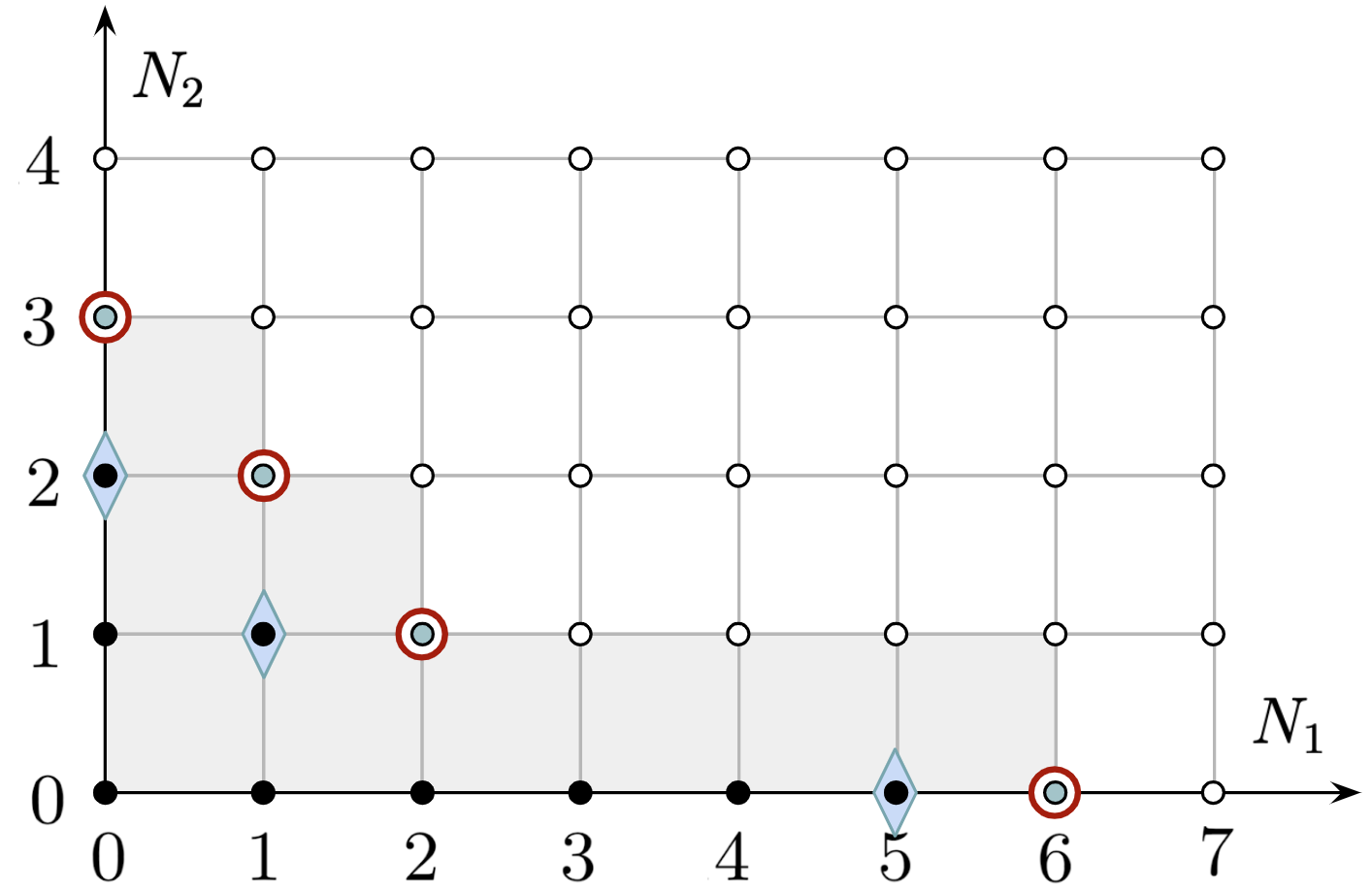}
\caption{The Young diagram of the partition $\lambda = (1 \leq 2 \leq 6)$ of $| \lambda|=9$ with addable corner set $\mathcal{A}_{\lambda} = \{(0,3), (1,2), (2,1), (6,0)\}$ (depicted above in circles), removable corner set $\mathcal{R}_{\lambda} = \{(0,2), (1,1), (5,0)\}$ (depicted above in diamonds), and outer corner set $\mathcal{R}_{\lambda}^+ = \{(1,3), (2,2), (6,1)\}$. Note that $\mathcal{R}_{\lambda} \subset \lambda$ but $\mathcal{R}_{\lambda}^+ \cap \lambda = \varnothing$.}\label{YoungDiagramFigure}
\end{figure}

\subsection{{The anisotropic content function}} \label{SUBSECAnisotropicContent} We now define {the following} function $[ \ \cdot \ ] \colon \BZ^2 \rightarrow \BC$.
\begin{Definition} For $\veps_1, \veps_2 \in \BC$, the anisotropic content $[s]$ of $s = (N_1, N_2) \in \BZ^2$ is defined by \begin{equation} \label{ContentDefinition} [s] = \veps_1 N_1 + \veps_2 N_2. \end{equation}
\end{Definition}
 For example, the anisotropic contents of the addable corner set $\mathcal{A}_{\lambda}$ of partition $\lambda {= (1 \leq 2 \leq 6)}$ in Figure~\ref{YoungDiagramFigure} are $3 \veps_2$, $\veps_1 + 2 \veps_2$, $2 \veps_1 + \veps_2$, and $6 \veps_1$. If $\veps_1 = \alpha$, $\veps_2 = -1$, these are $-3$, $\alpha -2$, $2 \alpha -1$, and $6 \alpha$, the {traditional anisotropic contents} of the {corresponding} $\alpha$-anisotropic Young diagram from~\cite{Ke4}. For simplicity, we may refer to $[s]$ as the content of $s$.

\subsection{Spectral theorem for the Nazarov--Sklyanin Lax operator} 
Nazarov--Sklyanin proved in~\cite{NaSk2} that the ingredients in Sections~\ref{SUBSECJackPolynomials}--\ref{SUBSECAnisotropicContent} emerge naturally in the spectral theory of the hierarchy $\mathcal{T}_{\ell} \colon F \rightarrow F$ defined by $\mathcal{T}_{\ell} = \pi_0 (\NSL_{\ebar,\hbar})^{\ell}$ provided one chooses $\veps_1, \veps_2 \in \BC$ so that
\begin{gather}
\label{EBarFormula} \ebar = \veps_1 + \veps_2, \\
\label{HBarFormula} \hbar = - \veps_1 \veps_2 .
\end{gather}

We will show that one can also discover partitions and Jack polynomials in the spectrum of~$\NSL$ itself.

\begin{Theorem}[main result] \label{TheoremMAINRESULT} Let $\NSL$ be the Nazarov--Sklyanin Lax operator in~\eqref{PresentationOfL}. Assume that $\ebar \in \BR$ and $\hbar>0$. Choose {non-zero} $\veps_1, \veps_2 \in \BC$ parametrizing $\ebar$, $\hbar$ by the formulas~\eqref{EBarFormula} and~\eqref{HBarFormula}.
\begin{enumerate}\itemsep=0pt
\item[$(1)$] The degree $n$ component $F[w]_n$ of the graded ring $F[w]$ in~\eqref{HDefinition} has a cyclic decomposition%
\begin{equation} \label{TheoremStatementDirectSum}
F[w]_n {\ = \ } \bigoplus_{| \lambda|=n} Z(j_{\lambda}, \NSL)
\end{equation}
as a direct sum of finite-dimensional $\NSL$-cyclic subspaces $Z(j_{\lambda},\NSL)$ indexed by parti\-tions $\lambda$ of size $n$ and generated by Jack polynomials $j_{\lambda} \in F$ defined in~\eqref{JackPolynomialPresentation}.
\item[$(2)$] The restrictions of $\NSL$ to the subspaces $Z_{\lambda} = Z(j_{\lambda}, \NSL)$ in~\eqref{TheoremStatementDirectSum} all have simple spectrum
\begin{equation}
\label{TheoremStatementSimpleSpectrum} \Spec \big (\NSL \big|_{Z_{\lambda}} \big ) = \big \{ [s] \colon s \in \mathcal{A}_{\lambda}\big \}
\end{equation}
given by the {anisotropic} contents~\eqref{ContentDefinition} of the addable corner set $\mathcal{A}_{\lambda}$ defined above in~\eqref{AddSetFormula}. {In particular, each cyclic subspace $Z_{\lambda}$ in~\eqref{TheoremStatementDirectSum} has dimension $\dim Z_{\lambda} = | \mathcal{A}_{\lambda}|$.}
\item[$(3)$] {As a distinguished vector space basis $\psi_{\lambda}^s$ for} the cyclic spaces $Z(j_{\lambda}, \NSL)$ in~\eqref{TheoremStatementDirectSum}, one can choose the unique polynomials $\psi_{\lambda}^s \in F[w]$ {in the variables $w, V_1, V_2, \ldots$} indexed {by addable corners} $s \in \mathcal{A}_{\lambda}$ which are both eigenfunctions of $\NSL$ with eigenvalues $[s]$ as in~\eqref{TheoremStatementSimpleSpectrum}, namely
\[
\NSL \psi_{\lambda}^s = [s] \psi_{\lambda}^s,
\]
{and which project to the Jack polynomial $j_{\lambda}$ for all $s \in \mathcal{A}_{\lambda}$ upon setting $w=0$, namely}
\begin{equation} \label{PsiNormalizationConditionTheoremStatement}
\pi_0 \psi_{\lambda}^s = j_{\lambda},
\end{equation}
 where $\pi_0\colon F[w] \rightarrow F$ is evaluation at $w=0$. In this normalization,
 \begin{equation} \label{JackSuperpositionTheoremStatement}
 j_{\lambda} = \sum_{s \in \mathcal{A}_{\lambda}} \tau_{\lambda}^s \psi_{\lambda}^s,
 \end{equation}
 where $\tau_{\lambda}^s>0$ are weights of the anisotropic transition measure \textnormal{\cite{Ke4}} defined explicitly by \begin{equation} \label{TauPresentation} \tau_{\lambda}^s = \frac{ \prod\limits_{r \in \mathcal{R}^+_{\lambda}} \ \ [s-r] }{ \prod\limits_{ t \in \mathcal{A}_{\lambda} \setminus \{s\}} [s - t] \ \ } \end{equation} for $\mathcal{R}_{\lambda}$ the removable corner set~\eqref{RemoveSetFormula} and $\mathcal{R}_{\lambda}^+ = \mathcal{R}_{\lambda}+(1,1)$ the set of outer corners.
\end{enumerate}
\end{Theorem}

Before we proceed, let us illustrate our Theorem~\ref{TheoremMAINRESULT} in degree $1$. For $\lambda=(1)$, $j_{1} = V_1$ and $Z_{1} = Z(j_1, \NSL)$ is two-dimensional with basis $\{V_1, w\}$. In this basis, $\NSL$ acts by the matrix
\begin{equation} \label{LRestrictedToZ1}
\NSL |_{Z_{1}} = \begin{bmatrix} 0 & 1 \\ \hbar & \ebar \end{bmatrix}.
\end{equation}
If $\ebar \in \BR$ and $\hbar>0$,~(\ref{LRestrictedToZ1}) has distinct real eigenvalues $\veps_1$ and $\veps_2$ determined up to permutation by the equations $\ebar = \veps_1 + \veps_2$ and $\hbar = - \veps_1 \veps_2$ characterizing the trace and determinant of~(\ref{LRestrictedToZ1}). In this way, the parametrization~(\ref{EBarFormula}) and~(\ref{HBarFormula}) above is visible already in degree $1$. Using~(\ref{ContentDefinition}), these simple eigenvalues $\veps_1 = [(1,0)]$, $\veps_2 = [(0,1)]$ are the contents of the two addable corners $\mathcal{A}_{{1}} = \{(0,1),(1,0)\}$ of {the unique partition} $\lambda={(1)}$ {of size $1$}. Moreover, as $\NSL |_{Z_1}$-eigenfunctions we can choose
\begin{gather*}
\psi_1^{(1,0)} = V_1 + \veps_1 w, \qquad \psi_1^{(0,1)} = V_1 + \veps_2 w
\end{gather*}
so that at $w=0$ one recovers $\psi_1^s \big |_{w=0} = j_{1}= V_1$ the Jack polynomial for both $s \in \mathcal{A}_{{1}} $.

\subsection{Organization of the paper} 
In Section~\ref{SECNazarovSklyanin}, we recall two results of Nazarov--Sklyanin~\cite{NaSk2}. In Appendix~\ref{APPENDIXKerovChapter1}, we collect standard results for cyclic spaces of self-adjoint operators in a self-contained appendix. In Section~\ref{SECProofSpectralTheorem}, we prove Theorem~\ref{TheoremMAINRESULT} using results from Section~\ref{SECNazarovSklyanin} and Appendix~\ref{APPENDIXKerovChapter1}. In Section~\ref{SECEigenvectors}, we {present} the polynomial eigenfunctions $\psi_{\lambda}^s$ of $\NSL$ {explicitly in degrees $| \lambda| \leq 3$}, {state an integrality conjecture for~$\psi_{\lambda}^s$, and prove a principal specialization formula for $\psi_{\lambda}^s$ which implies a special case of our integrality conjecture.} In Section~\ref{SECComments}, we discuss our proof and comment on related results in the literature.

\section{Review of two results of Nazarov--Sklyanin} \label{SECNazarovSklyanin}
 In this section, we recall two results for the hierarchy $\mathcal{T}_{\ell} = \pi_0 \NSL^{\ell}$ due to Nazarov--Sklyanin~\cite{NaSk2}.
\subsection{Definition of the Nazarov--Sklyanin hierarchy} 
Recall $F[w] = \BC[w, V_1, V_2, \ldots]$ in~(\ref{HDefinition}) and the Lax operator $\NSL$ in~(\ref{PresentationOfL}). Let $\pi_0$ be the projection to $F = \BC[V_1, V_2, \ldots]$ defined by setting $w=0$.

\begin{Definition} The Nazarov--Sklyanin hierarchy is the set of operators $\mathcal{T}_{\ell} = \pi_0 \NSL^{\ell} \colon F\rightarrow F$. \end{Definition}
Here $\mathcal{T}_{\ell}$ is the top-left entry of the $\ell^{\textnormal{th}}$ power of $\NSL$ in~(\ref{MatrixPresentationOfL}). Recall $V_{-k} = \hbar k \frac{\partial}{\partial V_k}$ in~(\ref{VNegativeKDefinition}). For~$\ell=3$, \begin{equation}
\label{DeformedCutAndJoin} \mathcal{T}_3 = \sum_{k_1, k_2 = 1}^{\infty} V_{-k_1 - k_2} V_{k_1} V_{k_2} + \sum_{k_1, k_2 = 1}^{\infty} V_{-k_1} V_{-k_2} V_{k_1 + k_2} + \ebar \sum_{k=1}^{\infty} k V_kV_{-k} \end{equation}
is the Hamiltonian of the quantum Benjamin--Ono equation on the torus discussed in~\cite[Section~1]{NaSk2}. When $\ebar = 0$,~(\ref{DeformedCutAndJoin}) is the well-known cut and join operator discussed, e.g., in~\cite[Section~1.2]{Dubrovin2014}.

\subsection{Commutativity and spectral theorem for the hierarchy} 

\begin{Theorem}[Nazarov--Sklyanin~\cite{NaSk2}] \label{NSTheorem}
Choose $\veps_1, \veps_2 \in \BC$ which parametrize $\ebar \in \BR$ and $\hbar >0$ by the formulas~\eqref{EBarFormula} and~\eqref{HBarFormula}. For $u \in \BC \setminus \BR$, consider the bounded transfer operators
\begin{equation}\label{NSTransferOperatorsPresentationIntro} \mathcal{T}(u) = \pi_0 (u - \NSL)^{-1} \colon \ F \rightarrow F, \end{equation}
 which encode all unbounded operators in the hierarchy $\mathcal{T}_{\ell} = \pi_0 \NSL^{\ell}$ as coefficients of $u^{-\ell-1}$.
\begin{itemize}\itemsep=0pt
\item[$(1^{\vee})$] The transfer operators~\eqref{NSTransferOperatorsPresentationIntro} commute in $F$ for distinct values of $u$ and are simultaneously diagonalized on Jack polynomials $j_{\lambda}$ in~\eqref{JackPolynomialPresentation} with eigenvalues $T_{\lambda}(u) \in \BC$ as in \begin{equation} \label{JacksFromNS} \mathcal{T}(u) \ j_{\lambda} = T_{\lambda}(u) \ j_{\lambda} .\end{equation}
\item[$(2^{\vee})$] The eigenvalues $T_{\lambda}(u) \in \BC$ in~\eqref{JacksFromNS} are determined by the {anisotropic} contents~\eqref{ContentDefinition} of the elements of the Young diagram $\lambda$ in~\eqref{DefinitionDLambda} and their shifts by the product formula
\begin{equation*} 
T_{\lambda}(u)= u^{-1} \cdot \prod_{s \in \lambda}\frac{(u-[s+(0,0)])(u-[s+(1,1)])}{(u-[s+(1,0)])(u-[s+(0,1)])}.
\end{equation*}

Equivalently, if $\mathcal{A}_{\lambda}$ and $\mathcal{R}_{\lambda}$ are the addable and removable corner sets in~\eqref{AddSetFormula} and~\eqref{RemoveSetFormula},%
\begin{equation}
\label{TEigenvalueAddRemoveStyle} T_{\lambda}(u) = \frac{\prod\limits_{r \in \mathcal{R}^+_{\lambda}} (u - [r]) }{ \prod\limits_{s \in \mathcal{A}_{\lambda} } (u - [s])} , \end{equation}
 where $\mathcal{R}_{\lambda}^+ = \mathcal{R}_{\lambda}+(1,1)$ is the set of outer corners.
 \end{itemize}
 \end{Theorem}

 Our presentation of the results of Nazarov--Sklyanin~\cite{NaSk2} above differs from what appears in~\cite{NaSk2}. For a proof that Theorem~\ref{NSTheorem} is equivalent to the original formulation in Nazarov--Sklyanin~\cite{NaSk2}, see~\cite[Section~8.2]{Moll2} and the discussion in~\cite[Section~4.3.3]{Moll5}. We label their two results $(1^{\vee})$ and $(2^{\vee})$ since we will now derive $(1)$ and $(2)$ in our Theorem~\ref{TheoremMAINRESULT} by proving that $(1^{\vee}) \land (2^{\vee}) \Rightarrow (1) \land (2)$. 

\section{Proof of main result} \label{SECProofSpectralTheorem}
 In this section we prove our Theorem~\ref{TheoremMAINRESULT}. Throughout, we assume that $j_{\lambda}$ in~(\ref{JackPolynomialPresentation}) and $\NSL$ in~(\ref{PresentationOfL}) have parameters $\ebar \in \BR$ and $\hbar>0$ satisfying $\ebar = \veps_1 + \veps_2$ as in~(\ref{EBarFormula}) and $\hbar = - \veps_1 \veps_2$ as in~(\ref{HBarFormula}).

\subsection{Jack--Lax cyclic spaces are finite-dimensional} 
We now recall the definition of certain cyclic spaces first considered in~\cite[Section~5.2.4]{Moll0}.
\begin{Definition} The Jack--Lax cyclic space $Z_{\lambda}$ is the $\NSL$-cyclic space generated by $j_{\lambda}$
\begin{equation} \label{JackLaxIntroductionInProof} Z_{\lambda} = Z(j_{\lambda}, \NSL), \end{equation} i.e., the subspace of $F[w]$ spanned by $\big\{ j_{\lambda}, \NSL j_{\lambda}, \NSL^2 j_{\lambda}, \NSL^3 j_{\lambda}, \ldots \big\}.$ \end{Definition}
A priori, we do not know $\dim Z_{\lambda}$ for $Z_{\lambda}$ in~(\ref{JackLaxIntroductionInProof}). To apply results from Appendix~\ref{APPENDIXKerovChapter1}, we need $\dim Z_{\lambda} < \infty$. The fact that Jack--Lax cyclic spaces $Z_{\lambda}$ are finite-dimensional follows from the next result.
\begin{Lemma} \label{JackLaxSpacesAreFiniteDim}The graded components $F[w]_n$ of $F[w]$ have dimension
\begin{equation} \label{HComponentDimensionFormula}
\dim F[w]_n= \sum_{l=0}^n \mathsf{p}(l),
\end{equation}
where $\mathsf{p}(l)$ is the number of partitions of $l$.
\end{Lemma}

\begin{proof} Write $\mu = \big(1^{d_1} 2^{d_2} \cdots\big)$ to denote a partition $\mu$ of size $| \mu | = \sum_{k=1}^{\infty} k d_k$ with $d_k$ parts of size $k$. For such $\mu$, let $V_{\mu} = V_1^{d_1} V_2^{d_2} \cdots$. The ring $F$ in~(\ref{H0Definition}) has graded components $F_n$ with $\dim F_n = \mathsf{p}(n)$ since $\{V_{\mu} \colon | \mu|=n\}$ is a basis of $F_n$. Similarly, $F[w]$ in~(\ref{HDefinition}) has graded components $F[w]_n$ with basis $\big\{V_{\mu} w^l \colon | \mu | + l = n\big\}$ so~(\ref{HComponentDimensionFormula}) holds.
\end{proof}

For $| \lambda|=n$, Jacks are homogeneous $j_{\lambda} \in F_n$ hence $j_{\lambda} \in F[w]_n$. Since $\NSL$ preserves degree,~$\NSL$~preserves $F[w]_n$, so $Z_{\lambda} = Z(j_{\lambda}, \NSL) \subset F[w]_n$ and~(\ref{HComponentDimensionFormula}) implies $\dim Z_{\lambda} < \infty$.

\subsection{Projections of Jack--Lax cyclic spaces} Next, we show that Jack polynomials diagonalizing the Nazarov--Sklyanin hierarchy in $(1^{\vee})$ of Theorem~\ref{NSTheorem} implies a remarkable property of $Z_{\lambda}$.
\begin{Lemma} For any partition $\lambda$, the projection of the Jack--Lax cyclic subspace in~\eqref{JackLaxIntroductionInProof} to $F \subset F[w]$ is one-dimensional and spanned by the Jack polynomial $j_{\lambda}$:
\begin{equation} \label{ProjectionJackLax} \pi_0 Z(j_{\lambda}, \NSL)= \BC j_{\lambda}.\end{equation}
\end{Lemma}

\begin{proof} For any $\eta \in Z(j_{\lambda}, \NSL)$, there is a polynomial $P$ so $\eta = P(\NSL) j_{\lambda}$. As a consequence,~$\pi_0 \eta$~is a finite linear combination of $\pi_0 \NSL^{\ell} j_{\lambda}$ indexed by $\ell \geq 0$. Expanding the resolvent $(u - \NSL)^{-1}= \sum_{\ell=0}^{\infty} u^{-\ell-1} \NSL^{\ell}$, part $(1^{\vee})$ of Theorem~\ref{NSTheorem} implies all $\pi_0 \NSL^{\ell} j_{\lambda} \in \BC j_{\lambda}$.
\end{proof}

We emphasize that this short proof does not require the explicit formula~(\ref{TEigenvalueAddRemoveStyle}) for the eigenvalues of the Nazarov--Sklyanin hierarchy in $(2^{\vee})$ of their Theorem~\ref{NSTheorem}, only commutativity in~$(1^{\vee})$.

\subsection[The Lax operator in Z\_lambda is self-adjoint for the extended Hall inner product]{The Lax operator in $\boldsymbol{Z_{\lambda}}$ is self-adjoint for the extended Hall inner product} 

Let $\langle \cdot, \cdot \rangle_{\hbar}$ be the inner product on $F[w]$ in which $V_{\mu} w^l$ for $\mu = \big(1^{d_1} 2^{d_2} \cdots\big)$ are orthogonal with
\begin{equation*} 
\big\| V_{\mu} w^l \big\|^2_{\hbar} = \prod_{k=1}^{\infty} (\hbar k)^{d_k} d_k!
\end{equation*}
By our discussion in Section~\ref{SUBSECJackPolynomials}, the restriction of $\langle \cdot , \cdot \rangle_{\hbar}$ to $F$ is the $\alpha$-Hall inner product from~\cite{Mac, Stanley}.

\begin{Lemma} \label{LRestrictionsSelfAdjoint} If $\ebar \in \BR$, $\hbar >0$, the restriction of $\NSL$ in~\eqref{PresentationOfL} to any $F[w]_n$ or $Z_{\lambda}$ is self-adjoint.
\end{Lemma}
\begin{proof} {Follows from the definition of $\NSL$ in~(\ref{PresentationOfL}), the fact that $V_{\pm k}^{\dagger} = V_{\mp k}$ in~(\ref{VNegativeKDefinition}) are mutual adjoints for $\langle \cdot, \cdot \rangle_{\hbar}$, and that the differential operator $w\partial_w$ is self-adjoint for $\langle \cdot, \cdot \rangle_{\hbar}$.}
\end{proof}

\subsection{Orthogonality of Jack--Lax cyclic spaces} We now use the projection formula in~(\ref{ProjectionJackLax}) and the orthogonality of Jack polynomials for the Hall inner product on $F$~\cite{Jack, Mac, Stanley} to prove that the Jack--Lax cyclic spaces themselves are orthogonal for the extended Hall inner product.

\begin{Lemma} \label{LEMMAJackLaxOrthogonality} If $\lambda$, $\gamma$ are distinct partitions of $n$, $Z_{\lambda}$ and $Z_{\gamma}$ are orthogonal in $(F[w]_n, \langle \cdot, \cdot \rangle_{\hbar})$.\end{Lemma}
\begin{proof} If $\eta \in Z_{\lambda} = Z(j_{\lambda}, \NSL)$, since $j_{\lambda}$ is cyclic, there is a polynomial $P$ so that $\eta = P(\NSL) j_{\lambda}$. Similarly, if $\varphi \in Z_{\gamma} = Z(j_{\gamma}, \NSL)$, there is a polynomial $Q$ so that $\varphi = Q(\NSL) j_{\gamma}$. Then
\[ \langle \varphi, \eta \rangle_{\hbar} = \langle Q(\NSL) j_{\gamma}, P(\NSL) j_{\lambda} \rangle_{\hbar} = \langle j_{\gamma}, Q(\NSL) P(\NSL) j_{\lambda} \rangle_{\hbar} = 0 \]
since $\NSL^{\dagger} = \NSL$ is self-adjoint in $F[w]_n$, $\pi_0 Q(\NSL) P(\NSL) j_{\lambda} \in \BC j_{\lambda}$ by~(\ref{ProjectionJackLax}), and $\langle j_{\gamma} , j_{\lambda} \rangle_{\hbar} = 0$.
\end{proof}

\subsection{Spectrum of the Lax operator in Jack--Lax cyclic spaces} We now show that spectrum of the hierarchy $\mathcal{T}_{\ell} = \pi_0 \NSL^{\ell}$ found by Nazarov--Sklyanin in $(2^{\vee})$ of their Theorem~\ref{NSTheorem} determines that of $\NSL$ in $Z_{\lambda}$. Precisely, the eigenvalues of $\NSL |_{Z_{\lambda}}$ are simple and given by the anisotropic contents $[s]$ of addable boxes $s \in \mathcal{A}_{\lambda}$.

\begin{proof}[Proof of part (2) of Theorem~\ref{TheoremMAINRESULT}] Since $\dim Z_{\lambda}<\infty$ by Lemma~\ref{JackLaxSpacesAreFiniteDim}, we may apply the general linear algebra results in the appendix Section~\ref{SUBSECrational} to the case $W = Z_{\lambda}$, $J = j_{\lambda}$, and $L = \NSL |_{Z_{\lambda}}$. From this perspective, formula~(\ref{TEigenvalueAddRemoveStyle}) in $(2^{\vee})$ of Nazarov--Sklyanin's Theorem~\ref{NSTheorem} gives an exact formula for the Titchmarsh--Weyl function $T(u)$ in~(\ref{TUDefinition}), so~(\ref{TAsRatioDeterminants}) implies~(\ref{TheoremStatementSimpleSpectrum}).
\end{proof}

\subsection[Cyclic decomposition of \protect{F[w]\_n} into Jack--Lax cyclic spaces Z\_lambda]{Cyclic decomposition of $\boldsymbol{F[w]_n}$ into Jack--Lax cyclic spaces $\boldsymbol{Z_{\lambda}}$}

\begin{proof}[Proof of part (1) of Theorem~\ref{TheoremMAINRESULT}] In any partition $\lambda$, the number of addable and removable corners always differ exactly by $1$:~(\ref{AddSetFormula}) and~(\ref{RemoveSetFormula}) satisfy
\begin{equation*} 
| \mathcal{A}_{\lambda} | = | \mathcal{R}_{\lambda}| +1.
\end{equation*}
 By induction in $n \geq 0$, it is straightforward to prove the combinatorial identities
\begin{equation} \label{Combo} \sum_{l=0}^n \mathsf{p}(l)= \sum_{| \lambda|=n} | \mathcal{A}_{\lambda}| = \sum_{| \nu| = n+1} | \mathcal{R}_{\nu}|,
\end{equation}
 where $\mathsf{p}(l)$ is the number of partitions of $l$. Indeed, $(\lambda, s) \mapsto (\lambda \cup \{s\} , s)$ defines a bijection
\[ \bigcup_{| \lambda|=n} \mathcal{A}_{\lambda} \rightarrow \bigcup_{| \nu|=n+1} \mathcal{R}_{\nu}. \]
By~(\ref{HComponentDimensionFormula}) and~(\ref{Combo}), $\dim F[w]_n = \sum_{| \lambda| = n}\! | \mathcal{A}_{\lambda}|$. By formula~(\ref{TheoremStatementSimpleSpectrum}) of part (2) of The\-o\-rem~\ref{TheoremMAINRESULT} proved above, since $\NSL |_{Z_{\lambda}}$ has simple spectrum indexed by $\mathcal{A}_{\lambda}$, $\dim Z_{\lambda} = | \mathcal{A}_{\lambda}|$. By the orthogonality in Lemma~\ref{LEMMAJackLaxOrthogonality}, $\bigoplus_{| \lambda|=n} Z_{\lambda} \subseteq F[w]_n$. However, since we've shown $\dim F[w]_n = \sum_{| \lambda|=n} \dim Z_{\lambda}$, we must have equality
\[ F[w]_n = \bigoplus_{| \lambda| = n} Z_{\lambda},
\]
which proves~(\ref{TheoremStatementDirectSum}).
\end{proof}

\subsection[Normalization of eigenfunctions of L]{Normalization of eigenfunctions of $\boldsymbol{\NSL}$}

At last, we can complete the proof of Theorem~\ref{TheoremMAINRESULT}.

\begin{proof}[Proof of part (3) of Theorem~\ref{TheoremMAINRESULT}] By Lemma~\ref{LRestrictionsSelfAdjoint}, we may apply the general linear algebra results in the appendices Appendices~\ref{SUBSECspectraltheorem} and~\ref{SUBSECkerovformulae} to the case $W = Z_{\lambda}$, $J = j_{\lambda}$, $L = \NSL |_{Z_{\lambda}}$, and $\langle \cdot, \cdot \rangle = \langle \cdot, \cdot \rangle_{\hbar}$. Throughout Appendix~\ref{APPENDIXKerovChapter1}, $\dim W = m+1$ and we index the $L$-eigenvector basis~$\psi^{(i)}$ by a superscript $i \in \{0,1,\ldots, m\}$. By the formula~(\ref{TheoremStatementSimpleSpectrum}) of part (2) of Theorem~\ref{TheoremMAINRESULT} proven above, $\dim Z_{\lambda} = | \mathcal{A}_{\lambda}|$, so we may instead index the $\NSL$-eigenvector basis $\psi_{\lambda}^s$ of $Z_{\lambda}$ by superscripts~${s \in \mathcal{A}_{\lambda}}$. In this way, our desired normalization condition~(\ref{PsiNormalizationConditionTheoremStatement}) on $\psi_{\lambda}^s ( w, V_1, V_2, \ldots)$ follows immediately from~(\ref{FirstJNormalizationCondition}). Likewise, the identities~(\ref{JackSuperpositionTheoremStatement}) and~(\ref{TauPresentation}) follow directly from~(\ref{JAsSuperposition}) and~(\ref{TauDefinition}) using the aforementioned identification of~(\ref{TEigenvalueAddRemoveStyle}) and~(\ref{TheoremStatementSimpleSpectrum}). This completes the proof of Theorem~\ref{TheoremMAINRESULT}.
\end{proof}
\begin{Remark} While we prove $(1^{\vee}) \land (2^{\vee}) \Rightarrow (1) \land (2)$, it is not hard to use the same general framework from Appendix~\ref{APPENDIXKerovChapter1} to prove the converse $(1) \land (2) \Rightarrow (1^{\vee}) \land (2^{\vee})$, so our Theorem~\ref{TheoremMAINRESULT} is in fact equivalent to results from~\cite{NaSk2}. That being said, neither the eigenvalues nor eigenfunctions of $\NSL$ are discussed explicitly in~\cite{NaSk2}. We give brief comments on the methods of proof in~\cite{NaSk2} in Section~\ref{SECComments}.
\end{Remark}

\section{On the eigenfunctions of the Nazarov--Sklyanin Lax operator} \label{SECEigenvectors}

 In this section, we consider the $\psi_{\lambda}^s \in F[w]$ of $\NSL$ normalized by the condition $\pi_0 \psi_{\lambda}^s= j_{\lambda}$ in part~(3) of our Theorem~\ref{TheoremMAINRESULT}. In Section~\ref{SUBSECMoreJack}, we present $\psi_{\lambda}^s$ explicitly in degrees $| \lambda| \leq 3$ . In Section~\ref{SUBSECSymmetry}, we discuss symmetries of $\psi_{\lambda}^s$ under $\veps_1 \leftrightarrow \veps_2$. In Section~\ref{SUBSECPsiIntegrality}, we state an integrality conjecture for $\psi_{\lambda}^s$. In Section~\ref{SUBSECPsiPrincipalSpecialization}, we verify a special case of this conjecture by using Stanley's principal specialization formula~\cite{Stanley} for $\psi_{\lambda}^s|_{w=0}=j_{\lambda}$ to derive a similar formula for $\psi_{\lambda}^s|_{w=1}$.
\subsection[Examples of eigenfunctions of L]{Examples of eigenfunctions of $\boldsymbol{\NSL}$} \label{SUBSECMoreJack}
Let $j_{\lambda} $ be the Jack polynomial in~(\ref{JackPolynomialPresentation}). Recall
\begin{gather}
j_{\varnothing} = 1 ,\nonumber\\
j_1 = V_1 ,\nonumber\\
j_{2} = V_1^2 + \veps_1 V_2 ,\label{Degree2Jack2}\\
j_{1,1} = V_1^2 + \veps_2 V_2 ,\label{Degree2Jack11}\\
j_3 = V_1^3 + 3 \veps_1 V_1 V_2 + 2 \veps_1^2 V_3 ,\nonumber\\
j_{1,2} = V_1^3 + (\veps_1 + \veps_2) V_1 V_2 + (\veps_1 \veps_2) V_3 ,\label{LookAtThisJack}\\
j_{1,1,1} = V_1^3 + 3 \veps_2 V_1 V_2 + 2 \veps_2^2 V_3.\nonumber
\end{gather}

 The eigenfunctions $\psi_{\lambda}^s$ of $\NSL$ are polynomials in $\BC[w,V_1,V_2,\ldots]$, normalized by $\pi_0 \psi_{\lambda}^s = j_{\lambda}$, and in degrees $ | \lambda | \leq 3$ are
\begin{gather}
\psi_{\varnothing}^{(0,0)} = j_{\varnothing}, \nonumber \\
 \psi_1^{(0,1)} = j_1 + \veps_2 w, \nonumber \\
 \psi_1^{(1,0)} = j_1 + \veps_1 w, \nonumber \\
 \psi_{2}^{(0,1)} = j_2 + \veps_2 V_1 w + \veps_1 \veps_2 w^2, \nonumber \\
 \psi_{2}^{(2,0)} = j_2 + 2 \veps_1 V_1 w + 2 \veps_1^2 w^2, \nonumber \\
 \psi_{1,1}^{(0,2)} = j_{1,1} + 2 \veps_2 V_1 w + 2 \veps_2^2 w^2, \nonumber \\
 \psi_{1,1}^{(1,0)} = j_{1,1} + \veps_1 V_1 w + \veps_1 \veps_2 w^2, \nonumber \\
\psi_3^{(0,1)} = j_3 + \veps_2 V_1^2 w + \veps_1 \veps_2 V_2 w + 2 \veps_1 \veps_2 V_1 w^2 + 2 \veps_1^2 \veps_2 w^3 , \nonumber \\
\psi_3^{(3,0)} = j_3 + 3 \veps_1 V_1^2 w + 3 \veps_1^2 V_2 w + 6 \veps_1^2 V_1 w^2 + 6 \veps_1^3 w^3 , \nonumber \\
 \psi_{1,2}^{(0,2)} = j_{1,2} + 2 \veps_2 V_1^2 w + \veps_1 \veps_2 V_2 w + \veps_2 ( \veps_1 + 2 \veps_2) V_1 w^2 + 2 \veps_1 \veps_2^2 w^3, \nonumber \\
\psi_{1,2}^{(1,1)} = j_{1,2} + (\veps_1 + \veps_2) V_1^2 w + \big(\veps_1^2 - \veps_1 \veps_2 + \veps_2^2\big)V_2 w +3 \veps_1 \veps_2 V_1 w^2 + \veps_1 \veps_2 (\veps_1 + \veps_2) w^3, \nonumber 
\\
\psi_{1,2}^{(2,0)} = j_{1,2} + 2 \veps_1 V_1^2 w +\veps_1 \veps_2 V_2 w + \veps_1 ( \veps_2 + 2 \veps_1) V_1 w^2 + 2 \veps_1^2 \veps_2 w^3, \label{Exemplar} \\
 \psi_{1,1,1}^{(0,3)} = j_{1,1,1} + 3 \veps_2 V_1^2 w + 3 \veps_2^2 V_2 w + 6 \veps_2^2 V_1 w^2 + 6 \veps_2^3 w^3, \nonumber \\
 \psi_{1,1,1}^{(1,0)} = j_{1,1,1} + \veps_1 V_1^2 w + \veps_1 \veps_2 V_2 w +2 \veps_1 \veps_2 V_1 w^2 + 2 \veps_1 \veps_2^2w^3. \nonumber 
\end{gather}

Readers more familiar with the $\alpha$ conventions in~\cite{Jack, Mac, Stanley} may set $\veps_1 = \alpha$ and $\veps_2 =-1$ so that $j_{\lambda} = J_{\lambda}$ in~(\ref{JackPolynomialPresentation}) {and $V_k = p_k$}. For example, if $J_{\lambda}$, $m_{\lambda}$, and $p_{\mu}$ are the Jack, monomial, and power sum symmetric functions in~\cite{Mac, Stanley}, $J_{1,1} = 2 m_{1,1} = p_1^2 - p_2$ and $J_{2} = (1 + \alpha) m_2 + 2 m_{1,1} = p_1^2 + \alpha p_2$ gives~(\ref{Degree2Jack2}) and~(\ref{Degree2Jack11}) by~(\ref{JackPolynomialPresentation}). At $\veps_1 = \alpha$, $\veps_2 = - 1$,~(\ref{LookAtThisJack}) is $J_{1,2} = p_1^3 + (\alpha-1) p_1 p_2 - \alpha p_3$ and so~(\ref{Exemplar}) is $\psi_{1,2}^{(2,0)} = p_1^3 + (\alpha-1) p_1 p_2 - \alpha p_3 + 2 \alpha p_1^2 w - \alpha p_2 w + \alpha (-1 + 2 \alpha) p_1 w^2 - 2 \alpha^2 w^3$.

\subsection[Symmetries of eigenfunctions of L]{Symmetries of eigenfunctions of $\boldsymbol{\NSL}$} \label{SUBSECSymmetry}
As is evident in the examples above, Jack polynomials $j_{\lambda}$ are invariant under simultaneous permutation $\veps_1 \leftrightarrow \veps_2$ and transposition $\lambda \leftrightarrow \lambda'$. This symmetry is easier to see in the modern conventions~\cite{LiQinWang, QinBOOK} in Nekrasov variables $\veps_1$, $\veps_2$~\cite{NekOk, Okounkov2018ICM}. Since the Lax operator $\NSL$ itself depends only on $\ebar = \veps_1 + \veps_2$ and $\hbar = - \veps_1 \veps_2$ which are both invariant under $\veps_1 \leftrightarrow \veps_2$, its eigenfunctions~$\psi_{\lambda}^s$ are invariant under simultaneous permutation $\veps_1 \leftrightarrow \veps_2$, transposition $\lambda \leftrightarrow \lambda'$, and reflection $s \leftrightarrow s'$, i.e., $(N_1, N_2) \leftrightarrow (N_2, N_1)$ in $ \BN^2$.

\subsection[Integrality conjecture for eigenfunctions of L]{Integrality conjecture for eigenfunctions of $\boldsymbol{\NSL}$} \label{SUBSECPsiIntegrality}
Although Jack polynomials $J_{\lambda}$ are referred to as the ``integral form'' Jack polynomials in~\cite{Mac, Stanley}, the fact that their coefficients in the power sum basis $p_{\mu}$ are polynomials in $\alpha$ with integer coefficients was not proven until~\cite{LapointeVinet}. To state this result for $j_{\lambda}$ as in~(\ref{JackPolynomialPresentation}), let $\chi_{\lambda, \mu}(\veps_1, \veps_2)$ denote the coefficient of $V_{\mu} = V_1^{d_1} V_2^{d_2} \cdots$ in
\begin{equation} \label{JackCharacterExpansion}
j_{\lambda}(V_1, V_2, \ldots ; \veps_1, \veps_2) = \sum_{\mu} \chi_{\lambda, \mu}(\veps_1, \veps_2) V_{\mu}.
\end{equation}
 These $\chi_{\lambda, \mu} ( \veps_1, \veps_2)$ are known as unnormalized {\it Jack characters}~\cite{DoFe4, Las1} since they are deformations of symmetric group characters which appear in the case $\ebar = \veps_1 + \veps_2 = 0$ of Schur polynomials.
\begin{Theorem}[Lapointe--Vinet~\cite{LapointeVinet}] \label{IntegralityJackResult} For any $\lambda$, $\mu$, in~\eqref{JackCharacterExpansion} one has $\chi_{\lambda, \mu}(\veps_1, \veps_2) \in \BZ[\veps_1, \veps_2]$.
\end{Theorem}
Consider the generalized Kostka coefficients $\kappa_{\lambda}^{\gamma}(\alpha)$ defined by the expansion $J_{\lambda} \!=\! \sum_{\gamma} \kappa_{\lambda}^{\gamma}(\alpha) \widetilde{m}_{\gamma}$ of the Jack polynomial $J_{\lambda}$ in the basis $\widetilde{m}_{\gamma}$ of augmented monomial symmetric functions~\cite{Mac}. Lapointe--Vinet~\cite{LapointeVinet} proved $\kappa_{\lambda}^{\gamma}(\alpha) \in \BZ[\alpha]$ which implies Theorem~\ref{IntegralityJackResult} since $\widetilde{m}_{\lambda}$ are polynomials in~$p_{\mu}$ with integer coefficients~\cite{Mac}. For the same reason, Theorem~\ref{IntegralityJackResult} is also a consequence of the proof of the Macdonald--Stanley conjecture $\kappa_{\lambda}^{\gamma}(\alpha) \in \BZ_{\geq 0} [\alpha]$ by Knop--Sahi~\cite{KnopSahi}.

We conjecture that an analog of Theorem~\ref{IntegralityJackResult} holds for the eigenfunctions ${\psi}_{\lambda}^s$ of the operator~$\NSL$. Let $\chi_{\lambda, \mu}^{s, l}(\veps_1, \veps_2)$ denote the coefficient of $V_{\mu} w^{l}$ in \begin{equation} \label{ExtendedJackCharacterExpansion} \psi_{\lambda}^s(w, V_1, V_2, \ldots ;\veps_1, \veps_2) = \sum_{l} \sum_{\mu} \chi_{\lambda, \mu}^{s,l} (\veps_1, \veps_2) V_{\mu} w^l .\end{equation}

\begin{Conjecture} \label{IntegralityConjecture} For any $\lambda$, $s$, $\mu$, and $l$ in~\eqref{ExtendedJackCharacterExpansion} one has $\chi_{\lambda, \mu}^{s, l} (\veps_1, \veps_2) \in \BZ[\veps_1, \veps_2].$\end{Conjecture}

 This integrality conjecture manifestly holds for the low degree examples of $\psi_{\lambda}^s$ presented in Section~\ref{SUBSECMoreJack}. For example, in the expansion of $\psi_{1,2}^{(2,0)}$ in~(\ref{Exemplar}), the coefficient of $V_1 w^2$ is $\chi_{(1,2), (1)}^{(2,0), 2} = \veps_1 ( \veps_2 + 2 \veps_1)$. Note that the operator $\NSL$ itself in~(\ref{PresentationOfL}) and~(\ref{MatrixPresentationOfL}) has coefficients in~$\BZ[\veps_1, \veps_2]$. In addition, by part (2) of our Theorem~\ref{TheoremMAINRESULT}, $\NSL$ has eigenvalues in $\BZ_{\geq 0} [ \veps_1, \veps_2]$ since $\Spec(\NSL) = \big\{ [s] \colon s \in \BN^2 \big\}$.

 We can verify Conjecture \ref{IntegralityConjecture} for the lowest and highest possible powers $l$ of $w$:
 \begin{enumerate}\itemsep=0pt
\item If $l=0$ so $V_{\mu} w^l = V_{\mu}$, the normalization condition $\pi_0 \psi_{\lambda}^s = j_{\lambda}$ implies $\chi_{\lambda, \mu}^{s,0}$ is the Jack character $\chi_{\lambda,\mu}$ in~(\ref{JackCharacterExpansion}), so $\chi_{\lambda, \mu}^{s,0} \in \BZ[\veps_1, \veps_2]$ holds as a consequence of Theorem~\ref{IntegralityJackResult}.
\item If $l = | \lambda|$ so $V_{\mu} w^l = w^{|\lambda|}$, Proposition~\ref{PsiPrincipalSpecializationFormula} below implies the content product formula\begin{equation} \label{PsiTop} \chi_{\lambda, 0}^{s,| \lambda|} = \prod_{t \in (\lambda \cup s)^{\times} } [t], \end{equation}
 where $(\lambda \cup s)^{\times} =\lambda \cup s \setminus \{(0,0)\}.$ Since $[t]$ are defined by~(\ref{ContentDefinition}), $ \chi_{\lambda, 0}^{s,| \lambda|} \in \BZ_{\geq 0} [ \veps_1, \veps_2]$.
\end{enumerate}

 In the next section, we will derive this formula~(\ref{PsiTop}) for the top degree coefficient of the eigenfunctions of $\NSL$. To illustrate this result, consider the partition $\lambda = (1,2)$ of size $3$ and its addable corner $(2,0) \in \mathcal{A}_{1,2}.$ While the elements of the Young diagram $\lambda = \{(0,0),(0,1),(1,0)\}$ have contents $0$, $\veps_2$, $\veps_1$, if we consider the addable corner $s = (2,0)$, then the product in~(\ref{PsiTop}) is over $(\lambda \cup s )^{\times}= \{(0,1),(1,0), (2,0)\}$ have contents $\veps_2$, $\veps_1$, $2 \veps_1$ whose product $2 \veps_1^2 \veps_2$ is the top degree coefficient of $w^3$ in the formula~(\ref{Exemplar}) for $\psi_{1,2}^{(2,0)}$.

\subsection[Principal specializations of eigenfunctions of L]{Principal specializations of eigenfunctions of $\boldsymbol{\NSL}$} \label{SUBSECPsiPrincipalSpecialization}
We now derive the closed formula~(\ref{PsiTop}) for the top coefficient of $w^{|\lambda|}$ in the expansion of $\psi_{\lambda}^s$ as a product of the contents in $(\lambda \cup s)^{\times}$. To do so, we derive \textit{two} principal specialization formulas for $\psi_{\lambda}^s$.
\begin{Proposition}[principal specializations] \label{PsiPrincipalSpecializationFormula} For $\lambda$ with $| \lambda|=n$, let $\psi_{\lambda}^s (w, V_1, V_2, \ldots, V_{n})$ be the polynomial eigenfunction of $\NSL$ with eigenvalue $[s]$ normalized by $\pi_0 \psi_{\lambda}^s = j_{\lambda}$. Fix $z \in \BC$. Then
\begin{gather} \label{PsiPrincipalSpecializationW0}
\psi_\lambda^s(0, z, z, \ldots, z) = \prod_{t \in \lambda} (z + [t]), \\
\label{PsiPrincipalSpecializationW1} \psi_\lambda^s(1, z, z, \ldots, z)= \prod_{t \in (\lambda \cup s)^{\times}} (z + [t])
\end{gather}
are two content product formulae for principal specializations at $w=0$ and $w=1$, respectively.
\end{Proposition}

\begin{proof} Since $\pi_0 \psi_{\lambda}^s = j_{\lambda}$, the $w=0$ principal specialization formula~(\ref{PsiPrincipalSpecializationW0}) is an immediate consequence of the well-known principal specialization result for $J_{\lambda}$ due to Stanley~\cite{Stanley}. To prove~(\ref{PsiPrincipalSpecializationW1}), let $\NSL\equiv\NSL_{\ebar,\hbar}$ be the Lax operator in~(\ref{PresentationOfL}), with the choices~(\ref{EBarFormula}) and~(\ref{HBarFormula}). The eigenvalue equation is \begin{equation} \label{Pi0Me} \NSL_{\ebar,\hbar} \psi_{\lambda}^s= [s] \psi_{\lambda}^s. \end{equation}
Taking $\pi_0$ of both sides of~(\ref{Pi0Me}), one can use $\pi_0 \psi_{\lambda}^s = j_{\lambda}$ on the right side, then observe that we can replace $\pi_0 \NSL_{\ebar,\hbar}$ by $\pi_0 \NSL_{0,0}$ on the left side. Indeed, by inspecting the terms in~(\ref{PresentationOfL}) with~$\hbar$ and~$\ebar$, one has $ \pi_0 \NSL_{\ebar,\hbar} \eta = \pi_0 \NSL_{0,0} \eta$ for any $\eta \in F[w]$. As a consequence,
\begin{equation} \label{ConsequencePi0}
\pi_0 \NSL_{0, 0} \psi_{\lambda}^s = [s] j_{\lambda}.
 \end{equation}
 Since $\pi_0 \NSL_{0,0} = \sum_{k=1}^{\infty} V_k \pi_0 w^{-k}$, if we expand our eigenfunction as
 \begin{equation} \label{PsiWPowerExpansion}
 \psi_{\lambda}^s = \sum_{l=0}^{n} ( \pi_{l} \psi_{\lambda}^s) w^l
 \end{equation}
 with each $ \pi_{l} \psi_{\lambda}^s \in F = \BC[V_1, V_2, \ldots]$ homogeneous of degree $n - l$,~(\ref{ConsequencePi0}) reads
\begin{equation} \label{Arrival} V_1 \pi_1 \psi_{\lambda}^s + V_2 \pi_2 \psi_{\lambda}^s + \cdots + V_n \pi_n \psi_{\lambda}^s = [s] j_{\lambda}(V_1, V_2, \ldots, V_n).
\end{equation}
On the other hand, the difference $\psi_{\lambda}^s(1) - \psi_{\lambda}^s(0)$ between $\psi_{\lambda}^s$ at $w=1$ and $w=0$ is \begin{equation} \label{W1EvalArrive} \pi_1 \psi_{\lambda}^s + \pi_2 \psi_{\lambda}^s + \cdots + \pi_n \psi_{\lambda}^s = \psi_{\lambda}^s (1, V_1, V_2, \ldots, V_n) - j_{\lambda} (V_1, V_2, \ldots, V_n).\end{equation}
Equation~(\ref{Arrival}) is an equality of polynomials in $V_k$ which are homogeneous of degree $n$, whereas equation~(\ref{W1EvalArrive}) is an equality of polynomials in $V_k$ which are not homogeneous. If we evaluate~$V_k = z$ for all $1 \leq k \leq n$ in both~(\ref{Arrival}) and~(\ref{W1EvalArrive}) and assume $z \neq 0$, the left sides of each differ by an overall factor of $z$, so we can conclude that
\begin{equation}
\label{PlopStanley} \psi_{\lambda}^s (1, z, z, \ldots, z) = \left ( 1 + \frac{[s]}{z} \right ) \cdot j_{\lambda}(z, z, \ldots, z) .
\end{equation}
Substituting Stanley's result~(\ref{PsiPrincipalSpecializationW0}) into~(\ref{PlopStanley}) yields~(\ref{PsiPrincipalSpecializationW1}).
\end{proof}

Using the notation $\pi_l \psi_{\lambda}^s$ for the coefficient of $w^l$ in $\psi_{\lambda}^s$ as in~(\ref{PsiWPowerExpansion}), setting $z=0$ in~(\ref{PsiPrincipalSpecializationW1}) yields the content product formula for the non-vanishing top degree coefficient $\pi_{n} \psi_{\lambda}^s$ presented in~(\ref{PsiTop}).

\section{Comments and comparison with previous results} \label{SECComments}

In this section, we discuss the larger context of our Theorem~\ref{TheoremMAINRESULT}. In Section~\ref{SUBSECNakamuraLaxOperator}, we comment on the classical Lax operator $\mathsf{L}$ of Nakamura~\cite{Nakamura1979} and Bock--Kruskal~\cite{BockKruskal1979} which served as the main inspiration for the construction of the Lax operator $\NSL$ by Nazarov--Sklyanin~\cite{NaSk2}. In Section~\ref{GKComparison}, we compare the appearance of Jack polynomials at $w=0$ in $\pi_0 \psi_{\lambda}^s = j_{\lambda}$ of our Theorem~\ref{TheoremMAINRESULT} to a~recent result of G\'{e}rard--Kappeler~\cite{GerardKappeler2019} for $\pi_0 \Phi_l$ of the classical $\mathsf{L}$ eigenfunctions $\Phi_l$. In Section~\ref{SUBSECNazarovSklyaninDiscussion}, we comment on developments since~\cite{NaSk2}. Finally, in Section~\ref{SUBSECGeometry} we compare the eigenvalue equation $\NSL \psi_{\lambda}^s = [s] \psi_{\lambda}^s$ in our Theorem~\ref{TheoremMAINRESULT} indexed by pairs $(\lambda, \lambda \cup s)$ of partitions which differ by $s$ to two other instances of this equation in the literature: (i) in type $A$ representation theory~\cite{MolevNazarovOlshanski, OkounkovVershik} at $\ebar = 0$ and (ii) in the equivariant cohomology of nested Hilbert schemes of points in $\BC^2$~\cite{LiQinWang, Nak0, Okounkov2018ICM, QinBOOK}.

\subsection{Comments on the Nakamura--Bock--Kruskal classical Lax operator}
\label{SUBSECNakamuraLaxOperator}
In~\cite[Sections~1 and~2]{NaSk2}, Nazarov--Sklyanin discuss how they thought to introduce their Lax operator $\NSL$ in~(\ref{PresentationOfL}) and~(\ref{MatrixPresentationOfL}) to the study of Jack polynomials $j_{\lambda}$.
They did so as a~natural consequence of two observations. On the one hand, Jack polynomials have been long known to be eigenfunctions of the Hamiltonian~(\ref{DeformedCutAndJoin}) of the quantum Benjamin--Ono equation on the torus~-- see~\cite[Section~2]{NaSk2} and~\cite{Moll2} and references therein. On the other hand, the classical Benjamin--Ono equation admits a Lax pair due to Nakamura~\cite{Nakamura1979} and Bock--Kruskal~\cite{BockKruskal1979}. When the spatial geometry is a torus, the Lax operator $\mathsf{L}$ from~\cite{BockKruskal1979, Nakamura1979} takes the form in~(\ref{PresentationOfClassicalL}) below. Assume $u \in L^2(\mathbb{T}, \BR)$ is a real-valued distribution with Fourier modes $u_k \in \BC$ satisfying $\sum_{k=1}^{\infty} |u_k|^2< \infty$, $u_{-k} = \overline{u}_k$, and $u_0 = 0$. Then with $\pi$ as in~(\ref{SzegoProjectionFormula}) and $\partial_w = \frac{\partial}{\partial w}$, \begin{equation} \label{PresentationOfClassicalL} \mathsf{L} = \veps_0 w \partial_w + \sum_{k=1}^{\infty} \overline{{u}_{k}} w^k + \sum_{k=1}^{\infty} u_k \pi w^{-k} \end{equation} is the Nakamura--Bock--Kruskal Lax operator for the classical Benjamin--Ono equation on the torus with dispersion coefficient $\veps_0 \in \BR$~\cite{BockKruskal1979, Nakamura1979}. This $\mathsf{L}$ is partially-defined on $\BC[w]$ and essentially self-adjoint with respect to the inner product $\langle \cdot, \cdot \rangle$ on $\BC[w]$ in which $w^l$ are an orthonormal basis. Just like~(\ref{PresentationOfL}), only the first terms involve the projection $\pi$ in (\ref{SzegoProjectionFormula}). However, unlike~(\ref{PresentationOfL}), the second infinite sum over $\overline{u_{k}} w^k$ applied to a polynomial in $\BC[w]$ will not yield a polynomial if infinitely-many $\overline{u_k} = u_{-k} \neq 0$. As in Section~\ref{SUBSECNazarovSklyaninLaxOperator}, we can write $\mathsf{L}$ as a $\BN \times \BN$ matrix \begin{equation} \label{MatrixPresentationOfClassicalL} \mathsf{L} =
 \begin{bmatrix} 0 \veps_0 & u_1 & u_2 & u_3 & u_4 & \cdots \\
\overline{u_1} & 1 \veps_0 & u_1 & u_2 & u_3 & \ddots \\
\overline{u_2} & \overline{u_1}& 2 \veps_0 & u_1 & u_2 & \ddots \\
\overline{u_3} & \overline{u_2} & \overline{u_1} & 3 \veps_0 & u_1 & \ddots \\
\overline{u_4} & \overline{u_3}& \overline{u_2} & \overline{u_1} & 4 \veps_0 & \ddots \\
\vdots & \ddots & \ddots & \ddots & \ddots & \ddots \end{bmatrix} .\end{equation}
With these two observations in mind, Nazarov--Sklyanin~\cite{NaSk2} realized that canonical quantization via~(\ref{VNegativeKDefinition}) \begin{equation} \label{CanonicalQuantization} (u_k, u_{-k}) \longrightarrow ( V_k , V_{-k} ) \end{equation}
performed directly in the matrix elements of $\mathsf{L}$ in~(\ref{MatrixPresentationOfClassicalL}) yields a well-defined $\NSL\colon F[w] \rightarrow F[w]$ in~(\ref{MatrixPresentationOfL}) whose powers are well defined \textit{without normal ordering}. Equivalently, Nazarov--Sklyanin realized that the classical field $u(x)$ can be directly replaced in $\mathsf{L}$ by the affine $\widehat{\mathfrak{gl}}_1$-current at level~$\hbar$ to get $\NSL$. For a discussion of~(\ref{CanonicalQuantization}) from the point of view of geometric quantization, see~\cite{Moll2}.

\subsection{Comparison to G\'{e}rard--Kappeler's action-angle coordinates} \label{GKComparison} This paper was inspired by recent spectral analysis~\cite{Gassot2021, Gerard2020Survey, GerardKappeler2019, GerardKappelerTopalov2020FA, GerardKappelerTopalov2022SURVEY, Moll1} of the Nakamura--Bock--Kruskal Lax operator $\mathsf{L}$ in~(\ref{MatrixPresentationOfClassicalL}). At $\veps_0 =0$,~(\ref{MatrixPresentationOfClassicalL}) is a Toeplitz operator whose spectrum has been studied for over a century~\cite{BottcherSilbermannIntro, SimonSzego}. At $\veps_0 \neq 0$, the spectrum of $\mathsf{L}$ is simple~\cite{GerardKappeler2019, Moll1}. In~\cite{Moll2, Moll1}, the second author proved that G\'{e}rard--Kappeler~\cite{GerardKappeler2019} independently found the classical limit of the quantum hierarchy of Nazarov--Sklaynin~\cite{NaSk2}. We can now make a second comparison to~\cite{GerardKappeler2019}. On the one hand, a main result of~\cite{GerardKappeler2019} is that the constant terms of the classical $\mathsf{L}$ eigenfunctions determine the action-angle coordinates of the classical Benjamin--Ono equation on the torus. On the other hand, in (3) of our Theorem~\ref{TheoremMAINRESULT}, we proved that the constant terms of the quantum~$\NSL$ eigenfunctions are Jack polynomials. We hope that the many structural results in G\'{e}rard--Kappeler~\cite{GerardKappeler2019} admit explicit quantizations which can shed further light on the objects in this paper.

\subsection{Comments on developments since the work of Nazarov--Sklyanin} \label{SUBSECNazarovSklyaninDiscussion} The study of $\mathcal{T}_{\ell}=\pi_0 \NSL^{\ell}$ and intricate proof of ($1^{\vee}$) and ($2^{\vee}$) in Theorem~\ref{NSTheorem} by Nazarov--Sklyanin~\cite{NaSk2} draws on their prior work~\cite{NaSk1} on the $N \rightarrow \infty$ limits of the Sekiguchi--Debiard operators $A_N^{(1)}, A_N^{(2)}, \ldots$ for Jack polynomials~\cite{Mac}. In~\cite{NaSk3, NaSk4}, Nazarov--Sklyanin generalized their results to the case of Macdonald polynomials. In~\cite{NaSk4}, they mention that their Theorem~\ref{NSTheorem} was independently discovered in a different form by Sergeev--Veselov~\cite{SergVes, SergVes2}. The relationship between the eigenfunctions $\psi_{\lambda}^s$ considered in this paper, the framework in~\cite{SergVes, SergVes2}, and the quantum Baker--Achiever function in~\cite{NaSk2} deserves further study.

\subsection{Comparison to spectral theorems in representation theory and geometry} \label{SUBSECGeometry}
 Our Theorem~\ref{TheoremMAINRESULT} is not the only appearance of contents of addable corners as eigenvalues of operators. In the case $\veps_1 = - \veps_2 = \veps$ so that $\ebar = 0$, $\hbar = \veps^2$, and $\alpha=1$, our eigenvalue equation degenerates to
 \begin{equation*} \NSL_{0, \hbar} \psi_{\lambda}^s = [s] \mathcal{\psi}_{\lambda}^s, \end{equation*}
where $\NSL_{0,\hbar}$ is a block Toeplitz operator, $\psi_{\lambda}^s$ are polynomials which recover Schur polynomials at $w=0$, and the anisotropic content is $\veps$ times the usual content. In this case, the content of a single addable corner is well known to describe the spectral theory at the heart of the representation theory of $\mathfrak{gl}_N$ and the symmetric group $S(n)$, see, e.g., discussions in Molev--Nazarov--Olshanski~\cite{MolevNazarovOlshanski} and Okounkov--Vershik~\cite{OkounkovVershik}.

For generic $\ebar = \veps_1 + \veps_2$, it is also well known that contents of addable corners arise in the equivariant cohomology of nested Hilbert schemes of points in the affine plane~\cite{LiQinWang, Nak0, Okounkov2018ICM, QinBOOK}. Let $X^{[n]}$ be the Hilbert scheme of $n$ points in $\BC^2$, i.e., all ideals $I \subset \BC [ x_1, x_2]$ so $\BC[x_1, x_2]/I$ is a~vector space of dimension $n$. Then
\[ X^{[n,n+1]} = \big\{ \big(I, \widetilde{I}\big) \in X^{[n]} \times X^{[n+1]} \colon I \subset \widetilde{I} \big\} \]
 is the nested Hilbert scheme of points in $\BC^2$. The tautological line bundle $L \rightarrow X^{[n,n+1]}$ has fibers
 \begin{equation} \label{TautologicalLineBundle}
 L |_{(I,\widetilde{I})} = \widetilde{I} / I.
 \end{equation}
 The action of the torus $\mathsf{T} = \BC^{\times} \times \BC^{\times}$ on $\BC^2$ induces an action on both $X^{[n]}$ and $X^{[n,n+1]}$. The $\mathsf{T}$-fixed points $I_{\lambda} \in X^{[n]}$ are monomial ideals indexed by partitions $\lambda$ with $| \lambda|=n$. The $\mathsf{T}$-fixed points $\eta_{\lambda}^s \in X^{[n,n+1]}$ are
\[ \eta_{\lambda}^s = (I_{\lambda}, I_{\lambda \cup s}) \]
indexed by pairs $(\lambda, s)$ with $| \lambda|=n$ and $s$ an addable corner in $ \mathcal{A}_{\lambda}$. For this reason, the classes $[I_{\lambda}]$ and $[\eta_{\lambda}^s]$ must be a basis in the $\mathsf{T}$-equivariant cohomology rings of $X^{[n]}$ and $X^{[n,n+1]}$, respectively. Moreover, if $s = (N_1, N_2) \in \BN^2$, the torus action with characters $\veps_1$, $\veps_2$ on the monomial $x_1^{N_1} x_2^{N_2} \in \BC[x_1, x_2]$ is determined by the content $[s] = \veps_1 N_1 + \veps_2 N_2$ in~(\ref{ContentDefinition}). As a~consequence, the operator $c_1(L) \cup -$ of cup product with the first Chern class of $L$ in~(\ref{TautologicalLineBundle}) must act diagonally in the basis $[\eta_{\lambda}^s]$ and satisfy the same eigenvalue equation as the Nazarov--Sklyanin Lax operator in our Theorem~\ref{TheoremMAINRESULT}. This observation suggests an extension of the isomorphism identifying Jack polynomials $j_{\lambda}$ with $[I_{\lambda}]$~\cite{LiQinWang, Nak0, QinBOOK} and $\mathcal{T}_3 = \pi_0 \NSL^{3}$ in~(\ref{DeformedCutAndJoin}) with the operator in Lehn~\cite{Lehn1999} to an identification of $\psi_{\lambda}^s \in F[w]$ introduced in this paper and $[\eta_{\lambda}^s]$.

\appendix

\section{Cyclic spaces of self-adjoint operators} \label{APPENDIXKerovChapter1}

In this appendix, we recall several standard results for cyclic spaces of self-adjoint operators following the treatment of Jacobi operators by Kerov~\cite[Section~6]{Ke1}. To streamline our proof in Section~\ref{SECProofSpectralTheorem}, we present these results without choosing an orthogonal basis in which our self-adjoint operators are tridiagonal.

\subsection[Cyclic spaces W, tilde W associated to an operator L with cyclic vector J]{Cyclic spaces $\boldsymbol{W}$, $\boldsymbol{\wtW}$ associated to an operator $\boldsymbol{L}$ with cyclic vector $\boldsymbol{J}$} \label{SUBSECcyclic}
Let $W$ be a $\BC$-vector space with $\dim_{\BC} W = m+1 < \infty$. Let $L\colon W \rightarrow W$ be a linear operator (not necessarily self-adjoint) and $J \in W$ a non-zero vector. Recall that $W$ is a {$L$-cyclic space generated by $J$} if the set of $m+1$ vectors $\big\{J, LJ, L^2 J, \ldots, L^{m} J\big\}$ is a basis for $W$. In this case, we say that $J$ is a {cyclic vector} for the operator $L$ and write $W = Z(J,L)$. In Proposition~\ref{PROPsecondcyclicspace} below, we recall a recipe which produces new cyclic spaces $\wtW=Z\big(\wtJ, \wtL\big)$ from a given cyclic space $W = Z(J,L)$. Let $\BC J$ denote the span of $J$. Choose any complementary subspace $\wtW \subset W$ of codimension $1$ so that
\begin{equation} \label{WDecomposition} W = \BC J + \wtW. \end{equation}
Let $\Pi_J$ and $\wtPi$ be the canonical projections onto $\BC J$ and $\wtW$, respectively, with $\ker \Pi_J = \wtW$ and $\ker \wtPi = \BC J$. For these canonical projections, one has $\textnormal{Id}_W = \Pi_J + \wtPi$.
Let $a \in \BC$ and $\wtJ \in \wtW$ be {uniquely} determined by the expansion of the vector $LJ$ {with respect to}~(\ref{WDecomposition}) as in
\begin{equation} \label{JTildeDefinition}
LJ = aJ + \wtJ,
\end{equation}
{so that $aJ = \Pi_J LJ$ and $\wtJ = \wtPi L J$.} Let $\wtL \colon \wtW \rightarrow \wtW$ be the linear operator with domain $\wtW$ defined by
\begin{equation} \label{LTildeDefinition}
\wtL = \wtPi L {|_{\wtW} },
\end{equation}
the restriction of $\wtPi L\colon W \rightarrow W$ to $\wtW$.

\begin{Proposition} \label{PROPsecondcyclicspace} Consider $\wtW$, $\wtJ$, $\wtL$ in~\eqref{WDecomposition}--\eqref{LTildeDefinition}. If $W = Z(J,L)$ is $L$-cyclic with cyclic vector~$J$, then the codimension $1$ subspace $\wtW$ is $\wtL$-cyclic with cyclic vector $\wtJ$, i.e., $\wtW=Z\big(\wtJ, \wtL\big)$.\end{Proposition}

\begin{proof} For any $\widetilde{\eta} \in \wtW$, we need to find a polynomial $\widetilde{P}$ so that
\[ \widetilde{\eta} = {\widetilde{P}}\big(\wtL\big) \wtJ . \]
Choose $\eta \in W$ {such that} $\wtPi \eta = \widetilde{\eta}$. Since $W = Z(J,L)$, there is a polynomial $P$ so that \begin{equation} \label{SlamThis} \eta =P(L) J. \end{equation}
 Apply $\wtPi$ to both sides of~(\ref{SlamThis}). {One can then} {use the identity $\textnormal{Id}_W = \Pi_J + \wtPi$,}~(\ref{JTildeDefinition}), and~(\ref{LTildeDefinition}) to prove that $\wtPi L^{\ell} J$ for $\ell = 0,1,2,\ldots$ are in the span of $\wtJ, \wtL \wtJ, \wtL^2 \wtJ, \ldots$ {by induction on $\ell$}, thus defining $\widetilde{P}$ from $P$.
\end{proof}

 \subsection[Two rational functions T, tilde T defined by L and J]{Two rational functions $\boldsymbol{T}$, $\boldsymbol{\wtT}$ defined by $\boldsymbol{L}$ and $\boldsymbol{J}$} \label{SUBSECrational} For $u \in \BC$ with both $u-L$ and $u - \wtL$ invertible, namely $u \not \in \big( \Spec(L) \cup \Spec(\wtL ) \big)$, let~$T(u)$ and~$\wtT(u)$ be the unique scalars for which
 \begin{gather} \label{TUDefinition}
 \Pi_J \frac{1}{u-L} J = T(u) J , \\
 \label{TTildeUDefinition} \Pi_J L \frac{1}{u - \wtL}\wtJ = \wtT(u) J.
 \end{gather}

 Formulas~(\ref{TUDefinition}) and~(\ref{TTildeUDefinition}) define two {meromorphic} functions $T$, $\wtT$ of $u \in \BC$ away from $\Spec(L)$, $\Spec\big(\wtL\big)$, respectively.

\begin{Proposition} 
Assume for simplicity that $a=0$ in~\eqref{JTildeDefinition} so that $LJ = \wtJ$. For all $u \in \BC$ at which $T(u) \neq 0$ and $\wtT(u)$ is well defined, the two functions $T$ and $\wtT$ in~\eqref{TUDefinition} and~\eqref{TTildeUDefinition} satisfy \begin{equation} \label{TViaHittingTimes} \frac{1}{T(u)} + \wtT(u) = u. \end{equation} {In particular, the left-hand side of~\eqref{TViaHittingTimes} is analytic in $u \in \BC$.}
\end{Proposition}
\begin{proof} The following argument is standard both in the theory of Jacobi matrices and in the study of lattice paths, see, e.g.,~\cite[Lemma 6.3.2]{Ke1}. However, this result is usually presented in the context of a choice of a tridiagonal matrix representation of $L$ and assuming that $L$ is self-adjoint, neither of which is necessary for the argument below. Since $a=0$, $u^{-2} \Pi_J L J = 0$. Expanding the resolvent using a geometric series,~(\ref{TUDefinition}) becomes
\begin{equation} \label{ExpandOleFaithful} T(u) J=u^{-1} J+ \sum_{\ell=2}^{\infty} u^{-\ell-1} \Pi_J L^{\ell} J. \end{equation}
For terms in~(\ref{ExpandOleFaithful}) with $\ell \geq 2$, use~(\ref{WDecomposition}) to insert $\ell-1$ copies of $\textnormal{Id}_W = \Pi_J + \wtPi$ between the $\ell$ copies of $L$ which appear in the expansion of $L^{\ell}$. This produces $2^{\ell-1}$ terms of the form $\Pi_J \Pi_1 L \Pi_2 L \Pi_3 \cdots L \Pi_{\ell-1} L J$ with each $\Pi_i \in \big\{\Pi_J, \wtPi\big\}$. Group these terms by the minimum value of $i$ so that $\Pi_i = \Pi_J$. Since $a=0$, we can assume $i \geq 2$. After relabeling indices, one checks~(\ref{ExpandOleFaithful}) becomes $T(u) = u^{-1} + u^{-1} \wtT(u) T(u)$ with $\wtT(u)$ as defined in~(\ref{TTildeUDefinition}).
 \end{proof}

The first function $T$ in~(\ref{TUDefinition}) is the ratio of the characteristic polynomials of $\wtL$ and $L$.

\begin{Proposition} \label{PROPcramer} The complex function $T(u)$ defined for $u \not\in \Spec(L) $ by~\eqref{TUDefinition} is \begin{equation} \label{TAsRatioDeterminants} T(u) = \frac{\det \big(u- \wtL\big)}{ \det (u - L)} \end{equation} a rational function with poles at $\Spec(L)$ and zeroes at $\Spec\big(\wtL\big)$. \end{Proposition}
\begin{proof} Apply Cramer's rule for $\Pi_J$ of solutions $\Phi(u)$ of the linear system $(u - L) \Phi = J$.
\end{proof}

 As a consequence of~(\ref{TViaHittingTimes}) and~(\ref{TAsRatioDeterminants}), the second function $\wtT$ is also a rational function.

\begin{Corollary} \label{COROLLARYhittingtimecramer} If $a=0$ in~\eqref{JTildeDefinition}, the function $\wtT(u)$ defined for $u \not\in \Spec\big(\wtL\big)$ by~\eqref{TTildeUDefinition} is
\begin{equation*} 
\wtT(u) = u - \frac{ \det (u-L)}{ \det \big(u-\wtL\big)}
\end{equation*}
a rational function with poles at $\Spec\big(\wtL\big)$.\end{Corollary}

\subsection[The case of self-adjoint L with cyclic vector J]{The case of self-adjoint $\boldsymbol{L}$ with cyclic vector $\boldsymbol{J}$} \label{SUBSECassumeselfadjoint} Introduce a Hermitian inner product $\langle \cdot, \cdot \rangle$ on the finite-dimensional complex vector space $W$. Adopt the convention $\langle \eta, \alpha \xi \rangle = \alpha \langle \eta, \xi \rangle$ that the inner product is $\BC$-linear in the right-most entry. We now revisit the results in Appendices~\ref{SUBSECcyclic} and~\ref{SUBSECrational} for $L$, $\wtL$ and $T$, $\wtT$ under the assumption that $L$ is a self-adjoint operator in the space $(W, \langle \cdot, \cdot \rangle)$.

 Given a non-zero vector $J \in W$, let $\wtW = (\BC J)^{\perp}$ be the orthogonal complement of $\BC J$. Let~$\Pi_J$ and $\wtPi = \Pi_J^{\perp}$ be the orthogonal projections so that~(\ref{WDecomposition}) is an orthogonal decomposition $W = \BC J \oplus \wtW = \BC J \oplus (\BC J)^{\perp}$.
By restriction, $\langle \cdot , \cdot \rangle$ defines an inner product on $\wtW =(\BC J)^{\perp}$. If $L$ is self-adjoint in $(W, \langle \cdot, \cdot \rangle)$, $\wtL = \wtPi L \wtPi$ in~(\ref{LTildeDefinition}) is self-adjoint in $(\wtW, \langle \cdot, \cdot \rangle)$ since $\wtPi$ is self-adjoint. Under these assumptions, we can give alternative formulas for the rational functions $T$ and $\wtT$: $\langle \cdot, \cdot \rangle$
\begin{gather} \label{TUinnerproductDefinition}
T(u) = \frac{ \big\langle J,(u-L)^{-1} J \big\rangle}{ \langle J, J \rangle} \\ \label{TTildeUinnerproductDefinition} \wtT(u) = \frac{\big\| \wtJ\big\|^2}{\|J\|^2} \cdot \frac{ \big\langle \wtJ, (u-\wtL)^{-1} \wtJ \big\rangle}{ \big\langle \wtJ, \wtJ \big\rangle},
\end{gather}
using~(\ref{TUDefinition}) and~(\ref{TTildeUDefinition}), and $\Pi_{J} \xi = \frac{\langle J, \xi \rangle }{ \langle J, J \rangle } J$. Here~(\ref{TTildeUinnerproductDefinition}) is $\wtT(u) = A T^{\vee}(u)$ where $A = \frac{\| \wtJ\|^2}{ \| J \|^2}>0$ and
\begin{equation} \label{CheckYourselfT}
T^{\vee}(u)= \frac{ \big\langle \wtJ, (u-\wtL)^{-1} \wtJ \big\rangle}{ \big\langle \wtJ, \wtJ \big\rangle} .
\end{equation}

 As we discussed in~\cite{Moll1}, we write~(\ref{TTildeUinnerproductDefinition}) without simplifying $\frac{ \| \wtJ\|^2}{\langle \wtJ, \wtJ \rangle}=1$ to clarify that $T^{\vee}(u)$ in~(\ref{CheckYourselfT}), not $\wtT(u)$ in~(\ref{TTildeUinnerproductDefinition}), is the Stieltjes transform of a {\it cotransition measure} in Kerov~\cite{Ke1, Ke4}. Note that $T$ and $\wtT$ are often called the \textit{Titchmarsh--Weyl functions} of $L$ and $\wtL$ -- see~\cite[Section~6]{Ke1}.

\subsection[Spectral theorem for self-adjoint operators L, tilde L with cyclic vectors J, tilde J]{Spectral theorem for self-adjoint operators $\boldsymbol{L}$, $\boldsymbol{\wtL}$ with cyclic vectors~$\boldsymbol{J}$,~$\boldsymbol{\wtJ}$} \label{SUBSECspectraltheorem}
Assume that $L$ and $\wtL$ are self-adjoint in $(W, \langle \cdot, \cdot \rangle)$ and $\wtW = (\BC J)^{\perp}$ with cyclic vectors $J$ and $\wtJ$ as in Appendix~\ref{SUBSECassumeselfadjoint}.
\begin{Theorem}[spectral theorem for cyclic spaces of self-adjoint operators] \label{GeneralSpectralTheorem} Let $L$, $\wtL$ be the self-adjoint operators in $W$, $\wtW$ of dimensions $m+1$, $m$ with cyclic vectors $J$, $\wtJ$ as in Appendix~\textnormal{\ref{SUBSECassumeselfadjoint}}. Then $L$, $\wtL$ have real eigenvalues {$\sigma^{(i)}$, $\widetilde{\sigma}^{(i)}$} which are simple and strictly interlacing
\[
\sigma^{(m)} < \widetilde{\sigma}^{(m)} < \sigma^{(m-1)} < \cdots < \widetilde{\sigma}^{(2)} < \sigma^{(1)} < \widetilde{\sigma}^{(1)} < \sigma^{(0)}.
\]
{Since the spectrum is simple, there exist bases $\big\{\psi^{(i)}\big\}_{i=0}^{m}$, $\big\{ \widetilde{\psi}^{(i)} \big\}_{i=1}^m$ of the cyclic spaces} $W$, $\wtW$, {which are} $L$, $\wtL$ eigenvectors
\begin{gather}
L \psi^{(i)} = \sigma^{(i)} \psi^{(i)}, \label{PreliminaryLEigenvector} \\
\wtL \widetilde{\psi}^{(i)} = \widetilde{\sigma}^{(i)} \widetilde{\psi}^{(i)} \label{PreliminaryLTildeEigenvector}
\end{gather}
{defined uniquely up to complex rescaling of each eigenvector. In fact, the normalizations of each eigenvector can be fixed by the canonical constraints}
\begin{gather} \label{FirstJNormalizationCondition}
\Pi_J \psi^{(i)} = J, \\ \label{SecondJNormalizationCondition} \Pi_J L \widetilde{\psi}^{(i)} = J .
\end{gather}
\end{Theorem}

\begin{proof} By the Cauchy interlacing theorem -- e.g.,~\cite[p. 242]{HornJohnsonBOOK} or~\cite[Section~6]{Ke1} -- it remains to prove~(\ref{FirstJNormalizationCondition}) and~(\ref{SecondJNormalizationCondition}). For any basis of $L$ eigenvectors $\psi^{(i)}$, since $\Pi_J$ projects onto $\BC J$,
\begin{equation} \label{GammaFormula}
\Pi_J \psi^{(i)}= \gamma^{(i)} J
\end{equation}
for some $\gamma^{(i)} \in \BC$ for all $i \in \{0,1,\ldots, m\}$. We first argue that $\gamma^{(i)}= 0$ is not possible in~(\ref{GammaFormula}), i.e., $\Pi_J \psi^{(i)} \neq 0$. Pairing~(\ref{GammaFormula}) with $J$, using the inner product and  the fact that the orthogonal projection $\Pi_J$ is self-adjoint, one arrives at the identity
\begin{equation} \label{LookAtThis}
\langle J, \psi^{(i)} \rangle = \gamma^{(i)} \| J\|^2.
\end{equation}
As a consequence, if we expand $J$ in the $\psi^s$ basis, formula~(\ref{LookAtThis}) determines the coefficients
\begin{equation} \label{ProdMeBrah}
J = \sum_{i = 0}^m \Bigg ( \gamma^{(i)} \frac{ \| J \|^2}{ \big\| \psi^{(i)}\big\|^2} \Bigg ) \psi^{(i)}.
\end{equation}
Since $J$ is cyclic for $L$ in $W$, its coefficients in the $\psi^{(i)}$ basis {are} non-zero for all $i \in \{0,1,\ldots, m\}$, thus~(\ref{ProdMeBrah}) guarantees $\gamma^{(i)} \neq 0$ for all $i \in \{0,1,\ldots, m\}$. Since the eigenvalues $\sigma^{(i)}$ of $L$ in $W$ are distinct, the $\psi^{(i)}$ are uniquely determined up to overall factors in $\BC^{\times}$, so we can choose these factors uniquely so $\gamma^{(i)} = 1$ in~(\ref{GammaFormula}), proving~(\ref{FirstJNormalizationCondition}). By Proposition~\ref{PROPsecondcyclicspace}, $\wtL$ is cyclic in $\wtW = (\BC J)^{\perp}$. Since $\wtL = \wtL^{\dagger}$, the same argument {above repeated in $\wtW$ instead of $W$, this time choosing \smash{$\widetilde{\gamma}^{(i)} =   \| J\|^2/ \| \wtJ \|^2 $} instead of $\gamma^{(i)} = 1$, guarantees} that if we consider $\wtJ \in \wtW$ in~(\ref{JTildeDefinition}), there is a unique basis $\widetilde{\psi}^{(i)}$ of $\wtL$-eigenvectors of $\wtW$ with
\begin{equation} \label{MiddleManFixa}
\Pi_{\wtJ} \widetilde{\psi}^{(i)} = \frac{\|J\|^2}{\big\| \wtJ\big\|^2} \wtJ .
\end{equation}
{Pairing both sides of~(\ref{MiddleManFixa}) with $\wtJ$ and using the fact that the orthogonal projection $\Pi_{\wtJ}$ to $\BC \wtJ$ is self-adjoint,~(\ref{MiddleManFixa}) implies
\begin{equation} \label{RevisedPaperEXHIBITA}
\big\langle \wtJ , \widetilde{\psi}^{(i)} \big\rangle = \| J \|^2 .
\end{equation} At the same time, since $\wtJ = \wtPi L J$ holds by~(\ref{JTildeDefinition}), $L^{\dagger} = L$, $\wtPi = \wtPi^{\dagger}$, and $\wtPi \widetilde{\psi}^{(i)} = \widetilde{\psi}^{(i)}$, we also have
\begin{equation} \label{RevisedPaperEXHIBITB}
\big\langle \wtJ, \widetilde{\psi}^{(i)} \big\rangle = \big\langle \wtPi L J, \widetilde{\psi}^{(i)} \big\rangle = \big\langle J , L \wtPi \widetilde{\psi}^{(i)} \big\rangle = \big\langle J , L \widetilde{\psi}^{(i)} \big\rangle.
\end{equation}
 Equating~(\ref{RevisedPaperEXHIBITA}) and~(\ref{RevisedPaperEXHIBITB}) implies} $\frac{ \langle {J} , L\widetilde{\psi}^{(i)} \rangle}{ \langle J, J \rangle } = 1$ which is equivalent to~(\ref{SecondJNormalizationCondition}).
\end{proof}

\subsection[Residues of T, tilde T and eigenvectors of self-adjoint L, tilde L]{Residues of $\boldsymbol{T}$, $\boldsymbol{\wtT}$ and eigenvectors of self-adjoint $\boldsymbol{L}$, $\boldsymbol{\wtL}$} \label{SUBSECkerovformulae}
Let $T$ and $\wtT$ be the rational functions associated to generic linear operators $L$ with cyclic vector~$J$ in Appendix~\ref{SUBSECrational}. Under the assumption that $L$ and $\wtL$ are self-adjoint as in Appendix~\ref{SUBSECassumeselfadjoint}, consider the residues
\begin{gather} \label{TauDefinition}
\tau^{(i)} = \textnormal{Res}_{u=\sigma^{(i)}} T(u), \\
\label{TauTildeDefinition} \widetilde{\tau}^{(i)} = \textnormal{Res}_{u=\widetilde{\sigma}^{(i)}} \wtT(u)
\end{gather}
of $T$ and $\wtT$ at the simple eigenvalues $\sigma^{(i)} \in \Spec(L)$ and $\widetilde{\sigma}^{(i)} \in \Spec\big(\wtL\big)$ from Appendix~\ref{SUBSECspectraltheorem}. We first recall how these residues appear in the calculation of the squared norms \smash{$\big\| \psi^{(i)} \big\|^2$}, \smash{$\big\| \widetilde{\psi}^{(i)} \big\|^2$} of the eigenvectors of $L$ and $\wtL$.

\begin{Proposition} 
Let $\psi^{(i)}$, $\widetilde{\psi}^{(i)}$ be the sets of eigenvectors of the self-adjoint operators $L$, $\wtL$ normalized with respect to the cyclic vector $J$ by the conditions \eqref{FirstJNormalizationCondition} and \eqref{SecondJNormalizationCondition}. Let $\tau^{(i)}$, $\widetilde{\tau}^{(i)}$ be the residues in \eqref{TauDefinition} and \eqref{TauTildeDefinition} of the rational functions $T(u)$, $\wtT(u)$ in \eqref{TUinnerproductDefinition} and~\eqref{TTildeUinnerproductDefinition}. Then
\begin{gather}
\tau^{(i)} \| \psi^{(i)}\|^2 = \|J\|^2, \label{PsiNormFormula}\\
\widetilde{\tau}^{(i)} \| \widetilde{\psi}^{(i)} \|^2 = \| J\|^2.\label{PsiTildeNormFormula}
\end{gather} In particular, {the residues $\tau^{(i)}$ and $\widetilde{\tau}^{(i)}$} in \eqref{TauDefinition} and \eqref{TauTildeDefinition} are {all} positive.
\end{Proposition}

\begin{proof} {For non-zero $\eta \in W$, let $\Pi_{\eta}$ be the orthogonal projection onto $\BC \eta$ defined by \smash{${\Pi_{\eta} \xi\! = \!\frac{\langle \eta, \xi \rangle}{ \langle \eta, \eta \rangle} \eta}$}. Since the orthogonal bases \smash{$\psi^{(i)}$}, \smash{$\widetilde{\psi}^{(i)}$} of $W$, $\wtW$ determine} resolutions of identities $\textnormal{Id}_W$, $\textnormal{Id}_{\wtW}${, we have}
\begin{gather} \label{ReturnOfTheMack} J = \sum_{i= 0}^{m} \Pi_{\psi^{(i)}} J = \sum_{i=0}^m \frac{\big\langle \psi^{(i)}, J \big\rangle }{ \big\| \psi^{(i)} \big\|^2 } \psi^{(i)}, \\
\label{ReturnOfTheMackTILDE}\wtJ = \sum_{i=1}^m \Pi_{\widetilde{\psi}^{(i)}} \wtJ = \sum_{i=1}^m \frac{\big\langle \widetilde{\psi}^{(i)}, \wtJ \big\rangle }{ \big\| \widetilde{\psi}^{(i)} \big\|^2 } \widetilde{\psi}^{(i)}.
\end{gather}

 {Substituting~(\ref{ReturnOfTheMack}) and~(\ref{ReturnOfTheMackTILDE}) into the inner product formulas~(\ref{TUinnerproductDefinition}) and~(\ref{TTildeUinnerproductDefinition}) for $T$, $\wtT$, one may use the eigenvalue equations~(\ref{PreliminaryLEigenvector}) and~(\ref{PreliminaryLTildeEigenvector}) and the orthogonality of eigenfunctions to get}
\begin{gather}
T(u) = \frac{1}{\| J\|^2} \sum_{i=0}^m \frac{1}{u - \sigma^{(i)} } \cdot \frac{ \big| \big\langle J, \psi^{(i)} \big\rangle \big|^2}{ \big\| \psi^{(i)} \big\|^2}, \label{TResolved} \\
\label{TTildeResolved} \wtT(u) = \frac{1}{\| J\|^2} \sum_{i=0}^m \frac{1}{u - \widetilde{\sigma}^{(i)} } \cdot \frac{ \big| \big\langle \wtJ, \widetilde{\psi}^{(i)} \big\rangle \big|^2}{\big \| \widetilde{\psi}^{(i)} \big\|^2}.
\end{gather}
 {Taking residues of~(\ref{TResolved}) and~(\ref{TTildeResolved}) at eigenvalues $\sigma^{(i)}$, $\widetilde{\sigma}^{(i)}$ gives
\begin{gather}
\tau^{(i)} = \frac{1}{\| J\|^2} \cdot \frac{ \big| \big\langle J, \psi^{(i)} \big\rangle \big|^2}{ \big\| \psi^{(i)} \big\|^2}, \label{TausResolved}\\
\widetilde{\tau}^{(i)} = \frac{1}{\| J\|^2} \cdot \frac{ \big| \big\langle \wtJ, \widetilde{\psi}^{(i)} \big\rangle \big|^2}{ \big\| \widetilde{\psi}^{(i)} \big\|^2}.\label{TauTildesResolved}
\end{gather}
 Finally, since the} normalization conditions~(\ref{FirstJNormalizationCondition}) and~(\ref{SecondJNormalizationCondition}) and $L = L^{\dagger}$ imply
\begin{gather}
\big\langle J, \psi^{(i)} \big\rangle = \| J\|^2, \label{JPsiInnerProduct}\\
\big\langle \wtJ, \widetilde{\psi}^{(i)} \big\rangle = \| J\|^2, \label{JTildePsiTildeInnerProduct}
\end{gather}
 {substituting~(\ref{JPsiInnerProduct}) and~(\ref{JTildePsiTildeInnerProduct}) into~(\ref{TausResolved}) and~(\ref{TauTildesResolved}) completes the proof.} \end{proof}

 In fact, these residues relate the eigenvectors and cyclic vectors themselves, not just their norms.
\begin{Proposition} \label{PROPResiduesInSuperpositions} In the self-adjoint case, the residues {$\tau^{(i)}$ and $\widetilde{\tau}^{(i)}$} of $T$ and $\wtT$ arise as coefficients in the expression of the cyclic vectors $J$ and $\wtJ$ as linear combinations of eigenvectors of~$L$ and~$\wtL$
\begin{gather}
J = \sum_{i=0}^m \tau^{(i)} \psi^{(i)}, \label{JAsSuperposition}\\
\wtJ = \sum_{i =1}^m \widetilde{\tau}^{(i)} \widetilde{\psi}^{(i)}\label{JTildeAsSuperposition}
\end{gather}
assuming that $\psi^{(i)}$ and $\widetilde{\psi}^{(i)}$ are normalized by the conditions~\eqref{FirstJNormalizationCondition} and \eqref{SecondJNormalizationCondition}. \end{Proposition}

\begin{proof} Formulas~(\ref{PsiNormFormula}) and~(\ref{JPsiInnerProduct}) imply $\tau^{(i)} = \frac{\langle \psi^{(i)}, J \rangle}{ \langle \psi^{(i)}, \psi^{(i)} \rangle}$. Hence~(\ref{JAsSuperposition}), since
\[
\Pi_{\psi^{{(i)}}} \xi = \frac{ \big\langle \psi^{(i)}, \xi \big\rangle}{\big\langle \psi^{(i)}, \psi^{(i)} \big\rangle} \psi^{(i)}
.\]
Similarly, combining formulas~(\ref{PsiTildeNormFormula}) and~(\ref{JTildePsiTildeInnerProduct}) implies~(\ref{JTildeAsSuperposition}).
\end{proof}

 We conclude with a `converse' of Proposition~\ref{PROPResiduesInSuperpositions}: the eigenvectors of $L$ and $\wtL$ themselves can be determined directly from the cyclic vectors $J$ and $\wtJ$ and the resolvents of $L$ and $\wtL$.
\begin{Proposition} 
Let $\psi^{(i)}$ and $\widetilde{\psi}^{(i)}$ be eigenvectors of self-adjoint operators $L$ and $\wtL$ normalized with respect to the cyclic vector $J$ by the conditions \eqref{FirstJNormalizationCondition} and \eqref{SecondJNormalizationCondition} with corresponding eigenvalues $\sigma^{(i)}$ and $\widetilde{\sigma}^{(i)}$. Then, assuming $a=0$ in \eqref{JTildeDefinition} so that $LJ = \wtJ$, the eigenvectors of $L$ and~$\wtL$ satisfy
\begin{gather}
\label{KerovFormulaForPsiTilde} \widetilde{\psi}^{(i)} = \frac{1}{L-\widetilde{\sigma}^{(i)}} J, \\
\label{KerovFormulaForPsi} \psi^{(i)} = J + \frac{1}{\sigma^{(i)} - \wtL} \wtJ .
\end{gather}
\end{Proposition}

\begin{proof} We first prove~(\ref{KerovFormulaForPsiTilde}). To begin, we observe that \begin{equation} \label{ThisExpressionLivesInWTilde} \Pi_J \frac{1}{ L - \widetilde{\sigma}} J = 0 \end{equation} since~(\ref{TUDefinition}) implies that the left-hand side of~(\ref{ThisExpressionLivesInWTilde}) is $-T(\widetilde{\sigma}^t) J$, but we know that $T(\widetilde{\sigma}^t)=0$ since \smash{$\widetilde{\sigma}^{(i)} \in \Spec\big(\wtL\big)$} is necessarily a zero of the rational function $T$ by Proposition~\ref{PROPcramer}. As a~consequence of~(\ref{ThisExpressionLivesInWTilde}), \smash{$\wtPi \frac{1}{ L - \widetilde{\sigma}^{(i)} } J = \frac{1}{ L - \widetilde{\sigma}^{(i)} } J $} which implies the eigenvalue relation
\[ \wtL \frac{1}{L - \widetilde{\sigma}^{(i)} } J = \wtPi L \wtPi \frac{1}{L - \widetilde{\sigma}^{(i)}} J = \wtPi \big(L- \widetilde{\sigma}^{(i)} + \widetilde{\sigma}^{(i)}\big) \frac{1}{ L - \widetilde{\sigma}^{(i)}} J = \widetilde{\sigma}^{(i)} \frac{1}{L - \widetilde{\sigma}^{(i)} } J . \]
Finally, since formula~(\ref{ThisExpressionLivesInWTilde}) also implies that \smash{$\Pi_J L \frac{1}{ L - \widetilde{\sigma}^{(i)}} J = \Pi_J J + \widetilde{\sigma}^{(i)} \Pi_J \frac{1}{ L - \widetilde{\sigma}^{(i)} } J = J$}, we have shown that \smash{$\frac{1}{L - \widetilde{\sigma}^{(i)}} J$} is an $\wtL$-eigenvector with eigenvalue $\widetilde{\sigma}^{(i)}$ satisfying~(\ref{SecondJNormalizationCondition}). Since these properties uniquely characterize $\widetilde{\psi}^{(i)}$, we have proven~(\ref{KerovFormulaForPsiTilde}). Next, we prove~(\ref{KerovFormulaForPsi}) assuming $a=0$ so that $LJ = \wtJ$ in~(\ref{JTildeDefinition}). To begin, the resolvent of $\wtL$ is well defined at $u=\sigma^{(i)}$ since in Theorem~\ref{GeneralSpectralTheorem} we have seen that $\Spec(L) \cap \Spec\big(\wtL\big)= \varnothing$. As a consequence, $\wtJ \in \wtW$ implies that $\frac{1}{\sigma^{(i)} - \wtL} \wtJ \in \wtW$, hence the vector on the right-hand side of~(\ref{KerovFormulaForPsi}) satisfies \[ \Pi_J \Bigg ( J + \frac{1}{\sigma^{(i)} - \wtL} \wtJ \Bigg ) = J \]
 the normalization condition~(\ref{FirstJNormalizationCondition}). To verify that this vector is indeed an $L$-eigenvector with eigenvalue $\sigma^{(i)}$, multiply by $L$ then calculate its projections onto $\BC J$ and $\wtW = ( \BC J)^{\perp}$. In the~$\BC J$ projection, we can simplify using our assumption $a=0$ in~(\ref{JTildeDefinition}), the definition~(\ref{TTildeUDefinition}) of $\wtT(u)$, and $\wtT\big(\sigma^{(i)}\big) = \sigma^{(i)}$ at the special value $u = \sigma^{(i)}$ due to Corollary~\ref{COROLLARYhittingtimecramer}:
 \begin{equation} \label{PiJOfLTimesRHS}
 \Pi_J L \Bigg ( J + \frac{1}{\sigma^{(i)} - \wtL} \wtJ \Bigg ) = \Pi_J L J + \Pi_J L \frac{1}{\sigma^{(i)} - \wtL} \wtJ = a J + T\big(\sigma^{(i)}\big) J = \sigma^s J.
 \end{equation} In the $\wtW = (\BC J)^{\perp}$
 projection, using $\frac{1}{\sigma^{(i)} - \wtL} \wtJ \in \wtW$, $\wtPi L \wtPi = \wtL$, and $\wtL = -\big({\sigma^{(i)}} - \wtL\big) + \sigma^{(i)}$,
\begin{gather}
\wtPi L \Bigg ( J + \frac{1}{\sigma^{(i)} - \wtL} \wtJ \Bigg ) = \wtPi L J + \wtL \frac{1}{\sigma^{(i)} - \wtL} \wtJ= \wtJ - \wtJ + \sigma^{(i)} \frac{1}{\sigma^{(i)} - \wtL} \wtJ \nonumber\\
\hphantom{\wtPi L \Bigg ( J + \frac{1}{\sigma^{(i)} - \wtL} \wtJ \Bigg )}{}= \sigma^{(i)} \frac{1}{\sigma^{(i)} - \wtL} \wtJ. \label{PiTildeOfLTimesRHS}
\end{gather}
Since~(\ref{PiJOfLTimesRHS}) and~(\ref{PiTildeOfLTimesRHS}) are $\sigma^{(i)}$ times the same projections of $J+ \frac{1}{ \sigma^{(i)} - \wtL} \wtJ$,~(\ref{KerovFormulaForPsi}) follows.
\end{proof}

\subsection*{Acknowledgements}

The authors would like to thank the referees for many helpful comments and suggestions. We would also like to express our sincere thanks to the staff at Darwin's Ltd. coffee and sandwich shop on Cambridge Street in Cambridge, MA for supporting our collaboration during the years 2014--2019.

\pdfbookmark[1]{References}{ref}
\LastPageEnding


\begin{thebibliography}{99}
\footnotesize\itemsep=0pt

\bibitem{AlexanderssonMicklerPAPER3}
Alexandersson P., Mickler R., New cases of the strong {S}tanley conjecture,
 {i}n preparation.

\bibitem{Bi1}
Biane P., Representations of symmetric groups and free probability,
 \href{https://doi.org/10.1006/aima.1998.1745}{\textit{Adv. Math.}} \textbf{138} (1998), 126--181.

\bibitem{BockKruskal1979}
Bock T.L., Kruskal M.D., A two-parameter {M}iura transformation of the
 {B}enjamin--{O}no equation, \href{https://doi.org/10.1016/0375-9601(79)90762-X}{\textit{Phys. Lett.~A}} \textbf{74} (1979),
 173--176.

\bibitem{BottcherSilbermannIntro}
B\"ottcher A., Silbermann B., Introduction to large truncated {T}oeplitz
 matrices, Universitext, \href{https://doi.org/10.1007/978-1-4612-1426-7}{Springer}, New York, 1999.

\bibitem{DoFe4}
Do{\l}\c{e}ga M., F\'eray V., On {K}erov polynomials for {J}ack characters,
 \href{https://doi.org/10.46298/dmtcs.2322}{\textit{Discrete Math. Theor. Computer Sci. Proc.}} \textbf{AS} (2013),
 539--550, \href{https://arxiv.org/abs/1201.1806}{arXiv:1201.1806}.

\bibitem{Dubrovin2014}
Dubrovin B., Symplectic field theory of a disk, quantum integrable systems, and
 {S}chur polynomials, \href{https://doi.org/10.1007/s00023-015-0449-2}{\textit{Ann. Henri Poincar\'e}} \textbf{17} (2016),
 1595--1613, \href{https://arxiv.org/abs/1407.5824}{arXiv:1407.5824}.

\bibitem{Gassot2021}
Gassot L., Zero-dispersion limit for the {B}enjamin--{O}no equation on the
 torus with bell shaped initial data, \href{https://doi.org/10.1007/s00220-023-04701-0}{\textit{Comm. Math. Phys.}} \textbf{401}
 (2023), 2793--2843, \href{https://arxiv.org/abs/2111.06800}{arXiv:2111.06800}.

\bibitem{Gerard2020Survey}
G\'erard P., A nonlinear {F}ourier transform for the {B}enjamin--{O}no equation
 on the torus and applications, \href{https://doi.org/10.5802/slsedp.138}{\textit{S\'emin. Laurent Schwartz, EDP Appl.}}
 \textbf{2019--2020} (2019--2020), 8, 19~pages.

\bibitem{GerardKappeler2019}
G\'erard P., Kappeler T., On the integrability of the {B}enjamin--{O}no
 equation on the torus, \href{https://doi.org/10.1002/cpa.21896}{\textit{Comm. Pure Appl. Math.}} \textbf{74} (2021),
 1685--1747, \href{https://arxiv.org/abs/1905.01849}{arXiv:1905.01849}.

\bibitem{GerardKappelerTopalov2020FA}
G\'erard P., Kappeler T., Topalov P., On the spectrum of the {L}ax operator of
 the {B}enjamin--{O}no equation on the torus, \href{https://doi.org/10.1016/j.jfa.2020.108762}{\textit{J.~Funct. Anal.}}
 \textbf{279} (2020), 108762, 75~pages, \href{https://arxiv.org/abs/2006.11864}{arXiv:2006.11864}.

\bibitem{GerardKappelerTopalov2022SURVEY}
G\'erard P., Kappeler T., Topalov P., On the {B}enjamin--{O}no equation on
 {$\mathbb{T}$} and its periodic and quasiperiodic solutions,
 \href{https://doi.org/10.4171/jst/398}{\textit{J.~Spectr. Theory}} \textbf{12} (2022), 169--193, \href{https://arxiv.org/abs/2103.09291}{arXiv:2103.09291}.

\bibitem{HoOb}
Hora A., Obata N., Quantum probability and spectral analysis of graphs,
 \textit{Theoret. and Math. Phys.}, \href{https://doi.org/10.1007/3-540-48863-4_1}{Springer}, Berlin, 2007.

\bibitem{HornJohnsonBOOK}
Horn R.A., Johnson C.R., Matrix analysis, 2nd ed., \href{https://doi.org/10.1017/CBO9780511810817}{Cambridge University Press},
 Cambridge, 2013.

\bibitem{Jack}
Jack H., A class of symmetric polynomials with a parameter, \href{https://doi.org/10.1017/S0080454100008517}{\textit{Proc. Roy.
 Soc. Edinburgh Sect.~A}} \textbf{69} (1970), 1--18.

\bibitem{Ke1}
Kerov S., Interlacing measures, in Kirillov's {S}eminar on {R}epresentation
 {T}heory, \textit{Amer. Math. Soc. Transl. Ser.}, Vol. 181, \href{https://doi.org/10.1090/trans2/181/02}{American
 Mathematical Society}, Providence, RI, 1998, 35--83.

\bibitem{Ke4}
Kerov S., Anisotropic {Y}oung diagrams and symmetric {J}ack functions,
 \href{https://doi.org/10.1007/BF02467066}{\textit{Funct. Anal. Appl.}} \textbf{34} (2000), 41--51,
 \href{https://arxiv.org/abs/math.CO/9712267}{arXiv:math.CO/9712267}.

\bibitem{KnopSahi}
Knop F., Sahi S., A recursion and a combinatorial formula for {J}ack
 polynomials, \href{https://doi.org/10.1007/s002220050134}{\textit{Invent. Math.}} \textbf{128} (1997), 9--22,
 \href{https://arxiv.org/abs/q-alg/9610016}{arXiv:q-alg/9610016}.

\bibitem{KvingeLicataMitchell}
Kvinge H., Licata A.M., Mitchell S., Khovanov's {H}eisenberg category, moments
 in free probability, and shifted symmetric functions, \href{https://doi.org/10.5802/alco.32}{\textit{Algebr. Comb.}}
 \textbf{2} (2019), 49--74, \href{https://arxiv.org/abs/1610.04571}{arXiv:1610.04571}.

\bibitem{LapointeVinet}
Lapointe L., Vinet L., A {R}odrigues formula for the {J}ack polynomials and the
 {M}acdonald--{S}tanley conjecture, \href{https://doi.org/10.1155/S1073792895000298}{\textit{Int. Math. Res. Not.}}
 \textbf{1995} (1995), 419--424, \href{https://arxiv.org/abs/q-alg/9509002}{arXiv:q-alg/9509002}.

\bibitem{Las1}
Lassalle M., Jack polynomials and free cumulants, \href{https://doi.org/10.1016/j.aim.2009.07.007}{\textit{Adv. Math.}}
 \textbf{222} (2009), 2227--2269, \href{https://arxiv.org/abs/0802.0448}{arXiv:0802.0448}.

\bibitem{Lehn1999}
Lehn M., Chern classes of tautological sheaves on {H}ilbert schemes of points
 on surfaces, \href{https://doi.org/10.1007/s002220050307}{\textit{Invent. Math.}} \textbf{136} (1999), 157--207,
 \href{https://arxiv.org/abs/math.AG/9803091}{arXiv:math.AG/9803091}.

\bibitem{LiQinWang}
Li W.-P., Qin Z., Wang W., The cohomology rings of {H}ilbert schemes via {J}ack
 polynomials, in Algebraic Structures and Moduli Spaces, \textit{CRM Proc.
 Lecture Notes}, Vol.~38, American Mathematical Society, Providence, RI, 2004,
 249--258, \href{https://arxiv.org/abs/math.AG/0411255}{arXiv:math.AG/0411255}.

\bibitem{Mac}
Macdonald I.G., Symmetric functions and {H}all polynomials, 2nd ed., \textit{Oxford
 Mathematical Monographs}, The Clarendon Press, Oxford University Press, New
 York, 1995.

\bibitem{MicklerPAPER2}
Mickler R., {J}ack {L}ittlewood--{R}ichardson coefficients and the
 {N}azarov--{S}klyanin {L}ax operator, \href{https://arxiv.org/abs/2306.11115}{arXiv:2306.11115}.

\bibitem{MolevNazarovOlshanski}
Molev A., Nazarov M., Ol'shanskii G., Yangians and classical {L}ie algebras,
 \href{https://doi.org/10.1070/RM1996v051n02ABEH002772}{\textit{Russian Math. Surveys}} \textbf{51} (1996), 205--282,
 \href{https://arxiv.org/abs/hep-th/9409025}{arXiv:hep-th/9409025}.

\bibitem{Moll0}
Moll A., Random partitions and the quantum {B}enjamin--{O}no hierarchy,
 \href{https://arxiv.org/abs/1508.03063}{arXiv:1508.03063}.

\bibitem{Moll2}
Moll A., Exact {B}ohr--{S}ommerfeld conditions for the quantum periodic
 {B}enjamin--{O}no equation, \href{https://doi.org/10.3842/SIGMA.2019.098}{\textit{SIGMA}} \textbf{15} (2019), 098, 27~pages,
 \href{https://arxiv.org/abs/1906.07926}{arXiv:1906.07926}.

\bibitem{Moll1}
Moll A., Finite gap conditions and small dispersion asymptotics for the
 classical periodic {B}enjamin--{O}no equation, \href{https://doi.org/10.1090/qam/1566}{\textit{Quart. Appl. Math.}}
 \textbf{78} (2020), 671--702, \href{https://arxiv.org/abs/1901.04089}{arXiv:1901.04089}.

\bibitem{Moll5}
Moll A., Gaussian asymptotics of {J}ack measures on partitions from weighted
 enumeration of ribbon paths, \href{https://doi.org/10.1093/imrn/rnab300}{\textit{Int. Math. Res. Not.}} \textbf{2023}
 (2023), 1801--1881, \href{https://arxiv.org/abs/2010.13258}{arXiv:2010.13258}.

\bibitem{Nak0}
Nakajima H., Lectures on {H}ilbert schemes of points on surfaces, \textit{Univ.
 Lecture Ser.}, Vol.~18, \href{https://doi.org/10.1090/ulect/018}{American Mathematical Society}, Providence, RI, 1999.

\bibitem{Nakamura1979}
Nakamura A., B\"acklund transform and conservation laws of the
 {B}enjamin--{O}no equation, \href{https://doi.org/10.1143/JPSJ.47.1335}{\textit{J.~Phys. Soc. Japan}} \textbf{47} (1979),
 1335--1340.

\bibitem{NaSk2}
Nazarov M., Sklyanin E., Integrable hierarchy of the quantum {B}enjamin--{O}no
 equation, \href{https://doi.org/10.3842/SIGMA.2013.078}{\textit{SIGMA}} \textbf{9} (2013), 078, 14~pages,
 \href{https://arxiv.org/abs/1309.6464}{arXiv:1309.6464}.

\bibitem{NaSk1}
Nazarov M., Sklyanin E., Sekiguchi--{D}ebiard operators at infinity,
 \href{https://doi.org/10.1007/s00220-013-1821-z}{\textit{Comm. Math. Phys.}} \textbf{324} (2013), 831--849, \href{https://arxiv.org/abs/1212.2781}{arXiv:1212.2781}.

\bibitem{NaSk3}
Nazarov M., Sklyanin E., Macdonald operators at infinity, \href{https://doi.org/10.1007/s10801-013-0477-2}{\textit{J.~Algebraic
 Combin.}} \textbf{40} (2014), 23--44, \href{https://arxiv.org/abs/1212.2960}{arXiv:1212.2960}.

\bibitem{NaSk4}
Nazarov M., Sklyanin E., Cherednik operators and {R}uijsenaars--{S}chneider
 model at infinity, \href{https://doi.org/10.1093/imrn/rnx176}{\textit{Int. Math. Res. Not.}} \textbf{2019} (2019),
 2266--2294, \href{https://arxiv.org/abs/1703.02794}{arXiv:1703.02794}.

\bibitem{NekOk}
Nekrasov N.A., Okounkov A., Seiberg--{W}itten theory and random partitions, in
 The Unity of Mathematics, \textit{Progr. Math.}, Vol. 244, \href{https://doi.org/10.1007/0-8176-4467-9_15}{Birkh\"auser},
 Boston, MA, 2006, 525--596, \href{https://arxiv.org/abs/hep-th/0306238}{arXiv:hep-th/0306238}.

\bibitem{Okounkov2018ICM}
Okounkov A., On the crossroads of enumerative geometry and geometric
 representation theory, in Proceedings of the {I}nternational {C}ongress of
 {M}athematicians~-- {R}io de {J}aneiro~2018. {V}ol.~{I}. {P}lenary lectures,
 \href{https://doi.org/10.1142/9789813272880_0030}{World Scientific Publishing}, Hackensack, NJ, 2018, 839--867,
 \href{https://arxiv.org/abs/1801.09818}{arXiv:1801.09818}.

\bibitem{OkounkovVershik}
Okounkov A., Vershik A., A new approach to representation theory of symmetric
 groups, \href{https://doi.org/10.1007/PL00001384}{\textit{Selecta Math.~(N.S.)}} \textbf{2} (1996), 581--605.

\bibitem{Ol1}
Olshanski G., Anisotropic {Y}oung diagrams and infinite-dimensional diffusion
 processes with the {J}ack parameter, \href{https://doi.org/10.1093/imrn/rnp168}{\textit{Int. Math. Res. Not.}}
 \textbf{2010} (2010), 1102--1166, \href{https://arxiv.org/abs/0902.3395}{arXiv:0902.3395}.

\bibitem{QinBOOK}
Qin Z., Hilbert schemes of points and infinite dimensional {L}ie algebras,
 \textit{Math. Surveys Monogr.}, Vol. 228, \href{https://doi.org/10.1090/surv/228}{American Mathematical Society},
 Providence, RI, 2018.

\bibitem{SergVes}
Sergeev A.N., Veselov A.P., Dunkl operators at infinity and {C}alogero--{M}oser
 systems, \href{https://doi.org/10.1093/imrn/rnv002}{\textit{Int. Math. Res. Not.}} \textbf{2015} (2015), 10959--10986,
 \href{https://arxiv.org/abs/1311.0853}{arXiv:1311.0853}.

\bibitem{SergVes2}
Sergeev A.N., Veselov A.P., Jack--{L}aurent symmetric functions, \href{https://doi.org/10.1112/plms/pdv023}{\textit{Proc.
 Lond. Math. Soc.}} \textbf{111} (2015), 63--92, \href{https://arxiv.org/abs/1310.2462}{arXiv:1310.2462}.

\bibitem{SimonSzego}
Simon B., Szeg\H{o}'s theorem and its descendants. {S}pectral theory for
 {$L^2$} perturbations of orthogonal polynomials, M.~B.~Porter Lectures,
 Princeton University Press, Princeton, NJ, 2011.

\bibitem{Stanley}
Stanley R.P., Some combinatorial properties of {J}ack symmetric functions,
 \href{https://doi.org/10.1016/0001-8708(89)90015-7}{\textit{Adv. Math.}} \textbf{77} (1989), 76--115.

\end{thebibliography}
\end{document}